\documentclass{article}
\usepackage{amssymb}
\usepackage{latexsym}
\usepackage{amsmath}
\usepackage{color}
\usepackage{CJK}
\usepackage{graphicx}
\usepackage{fancyhdr}
\usepackage{subfigure}

\newtheorem{theorem}{Theorem}[section]

\newtheorem{lemma}{Lemma}[section]
\newtheorem{proposition}{Proposition}[section]
\newtheorem{remark}{Remark}[section]

\numberwithin{equation}{section}
\newenvironment{prooff}{\medskip\par\noindent{\bf Proof.}\ }{\qquad
\raisebox{-0.5mm}{\rule{1.5mm}{4mm}}\vspace{6pt}}

\newcommand{\bbrn}{\mathbb{R}^N}
\newcommand{\bbr}{\mathbb{R}}

\newcommand{\bbn}{\mathbb{N}}

\begin{document}

\begin{CJK*}{GBK}{song}

\title{On Kirchhoff type equations with critical Sobolev exponent and Naimen's open problems}

\author{ Yisheng Huang$^{a}$\thanks{E-mail address: yishengh@suda.edu.cn (Yisheng Huang)}\quad Zeng Liu$^b$\thanks{Corresponding
author. E-mail address:  zliu@mail.usts.edu.cn; luckliuz@163.com
(Zeng Liu)}\quad
Yuanze Wu$^c$\thanks{E-mail address: {wuyz850306@cumt.edu.cn} (Yuanze Wu).}\\
\footnotesize$^a${\em Department\, of \,Mathematics,\, Soochow
\,University,\, Suzhou\, 215006 \,P.R. China}\\
\footnotesize$^b${\em  Department of Mathematics, Suzhou University of Science and Technology, Suzhou 215009, P.R. China}\\
\footnotesize$^c$ {\em  College\, of\, Sciences,\, China\,
University\, of\, Mining\, and\, Technology,\, \,Xuzhou \,221116\,
P.R. China } }
\date{}
\maketitle
\begin{abstract}
We study the following Brezis-Nirenberg problem of Kirchhoff type
$$
\left\{\aligned &-(a+b\int_{\Omega}|\nabla u|^2dx)\Delta u = \lambda|u|^{q-2}u + \delta |u|^{2}u, &\quad \text{in}\ \Omega, \\
&u=0,& \text{on}\ \partial\Omega,
\endaligned
\right.
$$
where $\Omega\subset \bbr^4$ is a bounded domain with the smooth
boundary $\partial\Omega$, $2\leq q<4$ and $a$, $b$, $\lambda$,
$\delta$ are positive parameters.  We obtain some new existence and
nonexistence results, depending on the values of the above
parameters, which improves some known results. The asymptotical
behaviors of the solutions are also considered in this paper.


\end{abstract}

\section{Introduction}\label{Sec01}
In this paper, we consider  the following Kirchhoff type problem
involving critical nonlinearity
\begin{equation}\label{eq01}
\left\{\aligned -(a+b\int_{\Omega}|\nabla u|^2dx)\Delta u  &= \lambda|u|^{q-2}u + \delta |u|^{2}u, &\quad \text{in}\ \Omega, \\
u&=0,& \text{on}\ \partial\Omega,
\endaligned
\right.
\end{equation}
where $\Omega\subset \bbr^4$ is a bounded domain with the smooth
boundary $\partial\Omega$, $2\leq q<4$, $a$, $b$, $\lambda$,
$\delta$ are positive parameters.

Eq.~\eqref{eq01} is related to the stationary of the Kirchhoff type
quasilinear hyperbolic equation
\begin{equation}\label{eq61}
u_{tt} - \Big(a+b\int_{\Omega}|\nabla u|^2dy\Big)\Delta u=f(x, u),
\quad \text{in}\ \Omega,
\end{equation}
that was first proposed by Kirchhoff in \cite{K83} describing the
transversal oscillations of a stretched string. For more details on
the physical and mathematical background of Eq.~\eqref{eq61}, we
refer the readers to the papers \cite{A12,K83,L78,MR03} and the
references therein.  Besides the physical motivation, such problems
are often referred to as being nonlocal because of the presence of
the term $\int_{\Omega}|\nabla u|^2dx$ which implies that the
equation in (1.1) is no longer a pointwise identity. This phenomenon
leads to some mathematical difficulties, which makes the study of
such a class of problems particularly interesting. Equations like
\eqref{eq01} has been studied extensively by using variational
methods recently, see
\cite{A12,A10,AF12,CKW11,F13,HZ12,HL14,LLS14,MR03,D14,D114,PZ06,WTXZ12,W15,WHL15,ZP06}
and the references therein.

In the study of Eq. \eqref{eq01}, our main concern is focused on the
nonlinearities with critical growth starting from the pioneering
paper \cite{BN83} by Brezis and Nirenberg, in which \eqref{eq01} is
studied in the case of $a=\delta=1$ and $b=0$.  Many mathematicians
have paid great attention to this type of critical problems for
their stimulating and challenging difficulties coming from the lack
of compactness of the Sobolev imbedding
$H^1_0(\Omega)\hookrightarrow L^{2^*}(\Omega)$, see e.g.
\cite{CFS84,CFP85,CP86,CSZ12,CW05,SZ10,W96,Z89} and the references
therein for the  existence and multiplicity results.  Recently the
Brezis-Nirenberg problem of Kirchhoff type is investigated in
\cite{A10,F13,D14,D114} and the references therein.

In order to describe better our results, we distinguish two cases of
$q=2$ and $2<q<4$ according to the range of $q$.

\subsection{The case of $q=2$ and $\lambda\geq a\lambda_1$, i.e. the indefinite case of Kirchhoff type Brezis-Nirenberg
problem}

In this case, Eq.~\eqref{eq01} becomes the
following Brezis-Nirenberg problem involving Kirchhoff type nonlocal
term
\begin{equation}\label{eq02}
\left\{\aligned-(a+b\int_{\Omega}|\nabla u|^2dx)\Delta u &= \lambda u + \delta |u|^{2}u, &\quad \text{in}\ \Omega, \\
u&=0,& \text{on}\ \partial\Omega,
\endaligned
\right.
\end{equation}
which has been proved to have a solution for $\lambda\in(0,
a\lambda_1)$ if and only if $0\leq b\mathcal{S}^2<\delta$ by  Naimen
in \cite{D14} recently, where $\lambda_1$ is the first eigenvalue of
$-\Delta$ on $\Omega$ and $\mathcal{S}$ is the best Sobolev
constant defined in \eqref{eq06}.   Then some natural questions are
arisen: what will happen if $\lambda\geq a\lambda_1$?  What about
the asymptotical behaviors of the solutions depending on the
parameters?  The first aim of this
paper is to try to answer  these questions.  Note that,
compared with the case of $\lambda< a\lambda_1$, the
case of $\lambda\geq a\lambda_1$ is more sophisticated because the
operator $-a\Delta-\lambda$ is indefinite.
 In order to deal with
this case, besides the typical difficulty caused by the lack of
compactness of the Sobolev embedding $H^1_0(\Omega)\hookrightarrow
L^4(\Omega)$, we will face the following three difficulties: the
first  is to prove the boundedness of the Palais-Smale sequence (the
$(PS)$ sequence in short); the second is the interaction between the
Kirchhoff type perturbation $\|u\|^4_{H^1_0(\Omega)}$ and the
critical nonlinearity $\int_{\Omega}|u|^4dx$, in particular, they
have the same exponent; the third is that the weak limit of the
bounded $(PS)$ sequence cannot be seen as the weak solution of
Eq.~\eqref{eq02} directly.

 According to
the relationship between $b\mathcal{S}^2$ and $\delta$, we will consider the following
two cases: $0<\delta\leq b\mathcal{S}^2$ and $0<b\mathcal{S}^2<
\delta$.  For the first case, in order to study the asymptotical
behaviors of the solutions, we give some existence results of solutions for the following
eigenvalue problem
\begin{equation}\label{eq49}
\left\{\aligned-(a+b\int_{\Omega}|\nabla u|^2dx)\Delta u &=\lambda u, &\quad \text{in}\ \Omega, \\
u&=0,& \text{on}\ \partial\Omega.
\endaligned
\right.
\end{equation}
\begin{theorem}\label{thm05}
Suppose that $a>0$, $b>0$ and $\lambda>0$.
\begin{itemize}
  \item [{\em(I)}] If $\lambda\in(0,a\lambda_1]$ then
  Eq.~\eqref{eq49} has no nontrivial solution.
  \item [{\em(II)}] If $\lambda\in (a\lambda_1, a\lambda_2]$ then
  Eq.~\eqref{eq49} has a unique positive ground state solution $u=\sqrt{\frac{\lambda-a\lambda_1}{b\lambda_1^2}}\varphi_{11}$.
  \item [{\em(III)}] If $\lambda\in (a\lambda_k, a\lambda_{k+1}]$ with $k\geq 2$
  then Eq.~\eqref{eq49} has a unique positive ground state solution
  $u=\sqrt{\frac{\lambda-a\lambda_1}{b\lambda_1^2}}\varphi_{11}$;
  moreover, Eq.~\eqref{eq49} has  sign-changing solutions
  $\overline{u}_j=\sqrt{\frac{\lambda-a\lambda_j}{b\lambda_j^2}}\psi_j$ for all $\psi_j\in M(\lambda_j)$
  satisfying
  $$
\psi_j=\sum_{i=1}^{i_j}c_{ij}\varphi_{ij}\ \text{with}\
\sum_{i=1}^{i_j}c^2_{ij}=1, \quad j=2,3,\cdots
  k,
  $$
  where $M(\lambda_j)$, $\varphi_{ij}$ and $i_j$ will be given in
 Section~\ref{Sec02}.
\end{itemize}
\end{theorem}
\begin{remark}
{\em As far as we know, there is no result on the eigenvalue
problem~\eqref{eq49} in literatures.  Comparing to the well-known results in the case
$b=0$, we see that Eq.~\eqref{eq49} also has infinitely many
sign-changing solutions provided $\lambda\in (a\lambda_k,
a\lambda_{k+1}]$ with $k\geq 2$.  Eq.~\eqref{eq49} can be seen as
the limited problem of \eqref{eq02} as $\delta\to 0^+$, which will help us to
study the asymptotical
behaviors of the solutions of \eqref{eq02}.}
\end{remark}

Now we give some existence and nonexistence results of solutions for  Eq.~\eqref{eq02} in the case of $0<\delta\leq
b\mathcal{S}^2$.
\begin{theorem}\label{thm01}
Suppose that $a>0$, $b>0$, $\lambda>0$ and $\delta>0$ with
$0<\delta\leq b\mathcal{S}^2$.
\begin{itemize}
  \item [{\em(1)}] If $\lambda\in(0,a\lambda_1]$ then
  Eq.~\eqref{eq02} has no nontrivial solution.
  \item [{\em(2)}] If $\lambda>a\lambda_1$ and $\delta<b\mathcal{S}^2$ then
  Eq.~\eqref{eq02} has one positive ground state solution.
  \item [{\em(3)}] Assume that $\lambda\in(a\lambda_k, a\lambda_{k+1}]$ with $k\geq 1$ and
  $\{\delta_n\}$ is a sequence of positive numbers satisfying $\delta_n\to
  0$ as $n\to \infty$.  Let $u_n$ be the positive ground state solution
  corresponding to $\delta_n$, then there exists $u_0\in
  H^1_0(\Omega)\setminus\{0\}$ such that $u_n\to u_0$ in
  $H^1_0(\Omega)$ as $n\to\infty$, up to a subsequence.  Moreover,
  $$
  u_0=\sqrt{\frac{\lambda-a\lambda_1}{b\lambda_1^2}}\varphi_{11}
  $$
  is a positive ground state solution of the problem of \eqref{eq49}.

\end{itemize}
\end{theorem}
\begin{remark}\label{rmk02}{\em
The conclusions of (3) of Theorem~\ref{thm01} gives a description of
asymptotical behavior of the positive ground state solution of
\eqref{eq02} as $\delta\to 0^+$, however, an interesting problem is
that, besides the positive ground state solution, is there  any
solution of \eqref{eq02} converging  to a sign-changing solution
of \eqref{eq49} as $\delta\to 0^+$ if $\delta<b\mathcal{S}^2$?}
\end{remark}

For the case of
$0<b\mathcal{S}^2<\delta$,  we have the following results.
\begin{theorem}\label{thm03}
Suppose that $a>0$, $b>0$, $\lambda>0$ and $\delta>0$ with
$\delta>b\mathcal{S}^2$.
\begin{itemize}
  \item [{\em(i)}] If
\begin{equation}\label{eq51}
a\lambda_1\leq \lambda\leq
\frac{a\lambda_1\delta}{\delta-b\mathcal{S}^2}
\end{equation}
then Eq.~\eqref{eq02} has at least one pair of
sign-changing  solutions under one of the following conditions:
\begin{itemize}
  \item [{\em(1)}] $\delta>0$, $b>0$ is sufficiently small;
  \item [{\em (2)}] $b>0$, $\delta>0$ is sufficiently large.
\end{itemize}
 \item [{\em(ii)}] For $a, \delta>0$ and let $\{b_n\}$ and $\{\overline{\lambda}_n\}$ be two
 sequences of positive numbers
satisfying $b_n\to 0$ and $\overline{\lambda}_n\to a\lambda_1$ as
$n\to \infty$.   Assume that  $u_n$ is the solution corresponding to
$b_n$ and $\overline{\lambda}_n$ obtained above, then there exists
$u_0\in H^1_0(\Omega)\setminus\{0\}$ such that $u_n\to u_0$ in
$H^1_0(\Omega)$ as $n\to\infty$, up to a subsequence, and $u_0$ is
a sign-changing solution of the following equation:
\begin{equation}\label{eq48}
\left\{\aligned -a\Delta u &= a\lambda_1 u + \delta |u|^{2}u, &\quad \text{in}\ \Omega, \\
u&=0,& \text{on}\ \partial\Omega.
\endaligned
\right.
\end{equation}
\end{itemize}
\end{theorem}
\begin{remark}\label{rmk01}{\em
\begin{itemize}
  \item[(i)]From our results, we see that when $\lambda\geq a\lambda_1$, the
cases of $b=0$ and $b>0$ are quite different.
  \item[(ii)]Theorems~\ref{thm01} and \ref{thm03} complement and generalize the
results of \cite[Theorem 1.1]{D14} and \cite[Theorem 0.1]{CFP85}.
The assumption \eqref{eq51} can be removed if $\lambda=a\lambda_1$.
Moreover, the asymptotical behaviors of the solutions are
respectively given in Theorem~\ref{thm01} as $\delta\to 0$ and in
Theorem~\ref{thm03} as $b\to 0$.
  \item[(iii)]According to \cite{CSZ12,CW05}, Eq.~\eqref{eq48} has multiple
nontrivial solutions if we take the place of $\lambda_1$ by
$\lambda$, moreover the existence of ground state solutions of Eq.~\eqref{eq48} also
obtained in \cite{SWW09,CSZ12}, however, both of the existence of
ground state solutions and the multiplicity of nontrivial solutions
of \eqref{eq02} are still unknown for $\lambda\geq a\lambda_1$ and
$\delta>b\mathcal{S}^2$.
\end{itemize}
}
\end{remark}

\subsection{The case of $2<q<4$ and Naimen's open question.}

The case of $2<q<4$ is more thorny and tough because the boundedness
of the (PS) sequences is hard to prove.  To overcome this
difficulty, by using the well-known monotonicity trick due to Struwe
(see also \cite{J99,S88}), Naimen~\cite{D14} obtained the following
theorem: \vspace{6pt}

\noindent{\bf Theorem~A.}\ (\cite[Theorem~1.6]{D14})\ {\em Suppose that
$2<q<4$. Let $b, \delta>0$ satisfy
$b\mathcal{S}^2<\delta<2b\mathcal{S}^2$ and let $\Omega\subset \bbr^4$
be strictly star shaped.  Furthermore, assume that one of the following
conditions $(C1)$, $(C2)$ and $(C3)$ holds:
\begin{itemize}
  \item [$(C1)$] $a>0$, $\lambda>0$ is small enough,
  \item [$(C2)$] $\lambda>0$, $a>0$ is large enough,
  \item [$(C3)$] $a>0$, $\lambda>0$ and $\delta/b>\mathcal{S}^2$ is
  sufficiently close to $\mathcal{S}^2$,
\end{itemize}
then Eq.~\eqref{eq01} has a nontrivial
solution.}\vspace{6pt}

\noindent {\bf Naimen's open question.}\quad On \cite[Page~1171]{D14},  Naimen asked whether conditions that $\mu<2b\mathcal{S}^2$, $\Omega\subset \bbr^4$ is strictly star-sharped and $(C1)$-$(C3)$ in Theorem~A are necessary to ensure the existence of the solutions of  Eq.~\eqref{eq01}.

 Another purpose of the present paper is to try to give an answer to Naimen's open question.  In order to do this, inspired by \cite{BJ07,CZ14, CZZ15,JS14}, we will construct a
bounded (PS) sequence and show that the critical level of the
functional is below the compactness threshold
$\frac{(a\mathcal{S})^2}{4(\delta-b\mathcal{S}^2)}$ by applying a
perturbation method, then we obtain the following result which gives a partial answer to Naimen's open question.

\begin{theorem}\label{thm04}
Suppose that $2<q<4$, $a>0$, $b>0$, $\lambda>0$ and $\delta>0$ with
$\delta>b\mathcal{S}^2$. Then there exists $b^*_1>0 $ such that for
each $b\in(0,b^*_1)$, Eq.~\eqref{eq01} has a solution if one of the
above conditions $(C1)$--$(C3)$ is satisfied.

\end{theorem}

It is worth observing that Theorem~\ref{thm04} is complementary to
the corresponding result of \cite[Theorem 1.1]{F13}, where $\lambda>0$ is assumed to be large enough.

Furthermore, according to Theorem~\ref{thm04}, a natural questions about Eq.\eqref{eq01} can be
arisen: what will happen if $0<\delta\leq
b\mathcal{S}^2$?
For this question, we have the following
nonexistence and multiplicity results.
\begin{theorem}\label{thm06}
Suppose that $2<q<4$, $a>0$, $b>0$, $\lambda>0$ and $\delta>0$ with
$\delta\leq b\mathcal{S}^2$.
\begin{itemize}
  \item [$(i)$] If $\delta< b\mathcal{S}^2$ and
\begin{equation}\label{eq1.5}
a\geq\frac{(4-q)|\Omega|^{\frac12}}{\mathcal{S}}\Big(\frac{\lambda}{2}\Big)^{\frac{2}{4-q}}\Big(\frac{q-2}{b\mathcal{S}^2-\delta}\Big)^{\frac{q-2}{4-q}}
\end{equation}
  then Eq.~\eqref{eq01} has no nontrivial solution;
  \item [$(ii)$] there exists $b^*_2>0$ dependent of $a$ and $\lambda$ such that for each $b\in (0, b^*_2)$,
  Eq.~\eqref{eq01} has at least two nontrivial solutions provided $\delta< b\mathcal{S}^2$
  and has at least one nontrivial solution provided $\delta= b\mathcal{S}^2$.
\end{itemize}
\end{theorem}

This paper is organized as follows.  In Section~\ref{Sec02} we give
some preliminaries.  In Section~\ref{Sec05} we give a proof of
Theorem~\ref{thm05}. Section~\ref{Sec03} is devoted to prove
Theorems~\ref{thm01} and \ref{thm03} by using Ekeland variational
principle and applying an idea of \cite{CFP85}, in order to overcome
the interaction between the nonlocal term and the critical
nonlinearity, we give some new estimates (see Lemma~\ref{lem04}). In
Section~\ref{Sec04} we prove Theorems~\ref{thm04} and \ref{thm06} by
using a perturbation method introduced in \cite{BJ07,JS14}.

\section{Preliminaries}\label{Sec02}

Throughout this paper, we will use the following common notations.

\noindent{\bf Notations:}
\begin{itemize}
  \item $L^p(\Omega)$ is the usual Lebesgue
space with the norm $ \|u\|_p=\big(\int_{\Omega}|u|^p
dx\big)^{1/p}$.
  \item $H^1_0(\Omega)$ and $\mathcal{D}^{1,2}(\bbr^4)$ are
the usual Sobolev space with the associated norms
$\|u\|=\big(\int_{\Omega}|\nabla u|^2dx\big)^{1/2}$ and
$\|u\|_{\mathcal{D}}=\big(\int_{\bbr^4}|\nabla u|^2dx\big)^{1/2}$
respectively.
  \item Denote by $0<\lambda_1<\lambda_{2}<\cdots<\lambda_{j}<\cdots, j\in \bbn$ the eigenvalues of $(-\Delta,
H^1_0(\Omega))$ and by $M(\lambda_j)$ the corresponding eigenspace
of $\lambda_j$ in $H^1_0(\Omega)$.
  \item $\{\varphi_{ij}\}$ is the orthogonal eigenfunctions corresponding
to $\lambda_j$ spanning $M(\lambda_j)$ and
$\int_{\Omega}|\varphi_{ij}|^2dx=1$, $i=1,2,\cdots i_j$, where
$i_j=\text{dim}M(\lambda_j)$.
  \item $\Omega\subset \bbrn$ is a bounded domain and denote $|\Omega|$
  by the Lebesgue measure of $\Omega$.
  \item $C, C_1, C_2, \ldots $ denote the positive (possibly
different) constants.
\end{itemize}

We first note that the corresponding functional of \eqref{eq01} $I:
H^1_0(\Omega)\to \mathbb{R}$ is given by
\begin{equation}\label{eq04}
I(u)=\frac{a}{2}\|u\|^2+\frac{b}{4}\|u\|^4-\frac{\lambda}{q}\|u\|^q_q-\frac{\delta}{4}\|u\|^4_4,
\ \forall\ u\in H^1_0(\Omega).
\end{equation}
Clearly, $I\in C^1(H^1_0(\Omega),\bbr)$ and
$$
I'(u)v=(a+b\|u\|^2)\int_{\Omega}\nabla u\nabla
vdx-\lambda\int_{\Omega}|u|^{q-2}uvdx-\delta\int_{\Omega}|u|^2uvdx,
\ \forall\ u, v\in H^1_0(\Omega).
$$

Recalling that the best Sobolev constant $\mathcal{S}$ of
$\mathcal{D}^{1,2}(\bbr^4)\hookrightarrow L^4(\bbr^4)$ defined as
\begin{equation}\label{eq06}
\mathcal{S}=\inf_{u\in\mathcal{D}^{1,2}(\bbr^4)\setminus\{0\}}\frac{\int_{\bbr^4}|\nabla
u|^2dx}{\big(\int_{\bbr^4}|u|^4dx\big)^{1/2}}
\end{equation}
can be achieved by the Talenti function
\begin{equation}\label{eq20}
U_{\varepsilon}(x):=\frac{(8\varepsilon)^{1/2}}{\varepsilon+|x-x_0|^2},
\end{equation}
which is the positive solution of the problem in whole space
$$
-\Delta U = U^3, \ \text{in}\ \bbr^4, \quad U(x)\to 0 \ \text{as}\
|x|\to\infty
$$
for some $\varepsilon>0$ and $x_0\in\bbr^4$.  As a matter of fact,
when $b\mathcal{S}^2<\delta$, the Talenti function multiplied by an
appropriate constant,
$$
W_{\varepsilon}:=\big(\frac{a}{\delta-b\mathcal{S}^2}\big)^{1/2}U_{\varepsilon}
$$
is nothing but a solution of the Kirchhoff type equation in whole
space
\begin{equation}\label{eq64}
-(a+b\int_{\bbr^4}|\nabla W|^2)\Delta W =\delta W^3, \ \text{in}\
\bbr^4, \quad W(x)\to 0, \ \text{as}\ |x|\to 0.
\end{equation}
Moreover, we can easily check that the functional energy of
$W_{\varepsilon}$ satisfies
$$
\frac{a}{2}\|W_{\varepsilon}\|^2_{\mathcal{D}}+\frac{b}{4}\|W_{\varepsilon}\|^4_{\mathcal{D}}
-\frac{\delta}{4}\int_{\bbr^4}|W_{\varepsilon}|^4dx=\frac{(a\mathcal{S}^2)}{4(\delta-b\mathcal{S}^2)}.
$$
However, it is easy to see that \eqref{eq64} has no nontrivial
solutions if $0<\delta\leq b\mathcal{S}^2$.  Indeed, if not, we
assume that $W\in \mathcal{D}^{1,2}(\bbr^4)\setminus\{0\}$ is a
solution of \eqref{eq64}, then we obtain
\begin{equation}\label{eq69}
0=a\|W\|^2_{\mathcal{D}}+b\|W\|^4_{\mathcal{D}}-\delta\int_{\bbr^4}|W|^4dx\geq
a\|W\|^2_{\mathcal{D}}+(b-\frac{\delta}{\mathcal{S}^2})\|W\|^4_{\mathcal{D}}\geq
a\|W\|^2_{\mathcal{D}}>0,
\end{equation}
this is a contradiction.

The following result can be obtained using the mathematical
induction.
\begin{lemma}\label{lem05}
Let $x_k$, $y_k\in \bbr^+\ (k=1, 2, \cdots n)$ be positive
constants, then it holds
$$
\frac{\sum_{i=1}^{n}x_i}{\sum_{i=1}^{n}y_i}\leq
\max\Big\{\frac{x_k}{y_k}, k=1,2,\cdots n\Big\}.
$$
\end{lemma}

Next, we list a global compactness result which gives the
description of $(PS)_c$ sequence (cf. \cite{D14}).
\begin{proposition}\label{pro01}
Let $c\in\bbr$ and $\{u_n\}\subset H^1_0(\Omega)\subset
\mathcal{D}^{1,2}(\bbr^4)$ be a bounded $(PS)_c$ for $I$, that is,
$I(u_n)\to c$, $I'(u_n)\to 0$ in $H^{-1}(\Omega)$ and $\{u_n\}$ is
bounded.  Then $\{u_n\}$ has a subsequence which strongly converges
in $H^1_0(\Omega)$, or otherwise, there exist a function $u_0\in
H^1_0(\Omega)$ which is a weak convergence of $\{u_n\}$, a number
$k\in \bbn$ and further, for every $i\in\{1, 2, \cdots k\}$, a
sequence of value $\{R^i_n\}_{n\in\bbn}\subset \bbr^+$, points
$\{x_i\}_{n\in\bbn}\subset\overline{\Omega}$ and a function $v_i\in
\mathcal{D}^{1,2}(\bbr^4)$ satisfying
\begin{equation}\label{eq07}
-\Big\{a+b\Big(\|u_0\|^2+\sum_{j=1}^k\|v_j\|^2_{\mathcal{D}}\Big)\Big\}\Delta
u_0 = \lambda u_0^q+\delta u_0^3, \quad \text{in}\ \Omega,
\end{equation}
\begin{equation}\label{eq08}
-\Big\{a+b\Big(\|u_0\|^2+\sum_{j=1}^k\|v_j\|^2_{\mathcal{D}}\Big)\Big\}\Delta
v_i = \delta v_i^3, \quad \text{in}\ \bbr^4,
\end{equation}
such that up to subsequences, there hold $R^i_ndist(x^i_n,
\partial\Omega)\to\infty$,
$$
\Big\|u_n-u_0-\sum_{i=1}^kR^i_nv_i(R^i_n(\cdot-x^i_n))\Big\|_{\mathcal{D}}=o(1),
$$
$$
\|u_n\|^2=\|u_0\|^2+\sum_{i=1}^k\|v_i\|^2_{\mathcal{D}}+o(1)
$$
and
$$
I(u_n)=\widetilde{I}(u_0)+\sum_{i=1}^k\widetilde{I}^{\infty}(v_i)+o(1)
$$
where $o(1)\to 0$ as $n\to\infty$ and we define
\begin{equation}\label{eq09}
\widetilde{I}(u_0):=\Big\{\frac{a}{2}+\frac{b}{4}\Big(\|u_0\|^2+\sum_{j=1}^k\|v_j\|^2_{\mathcal{D}}\Big)\Big\}\|u_0\|^2
-\frac{\lambda}{q}\int_{\Omega}|u_0|^qdx-\frac{\delta}{4}\int_{\Omega}|u_0|^4dx,
\end{equation}
\begin{equation}\label{eq10}
\widetilde{I}^{\infty}(u_0):=\Big\{\frac{a}{2}+\frac{b}{4}\Big(\|u_0\|^2+\sum_{j=1}^k\|v_j\|^2_{\mathcal{D}}\Big)\Big\}\|v_i\|^2_{\mathcal{D}}
-\frac{\delta}{4}\int_{\Omega}|v_i|^4dx.
\end{equation}
\end{proposition}

Finally, we  recall a critical point theorem (cf. \cite[Theorem
2.4]{BBF83}) which is a variant of some results contained in
\cite{AR73}.
\begin{theorem}\label{thm02}
Let $H$ be a real Hilbert space and $I\in C^1(H, \bbr)$ be a
functional satisfying the following assumptions:
\begin{itemize}
  \item [$(a_1)$] $I(u)=I(-u)$, $I(0)=0$ for any $u\in H$;
  \item [$(a_2)$] there exists $\beta>0$ such that $I$ satisfies
  $(PS)$ condition in $(0, \beta)$;
  \item [$(a_3)$] there exist two closed subspace $V$, $W\subset H$ and
  positive constants $\rho$, $\eta$ such that
  \begin{itemize}
    \item [$(a_{31})$] $I(u)<\beta$ for any $u\in W$;
    \item [$(a_{32})$] $I(u)\geq \eta$ for any $u\in V$ with $\|u\|=\rho$;
    \item [$(a_{33})$] codim$V<\infty$.
  \end{itemize}
\end{itemize}
Then there exist at least $m$ pair of critical points with the
corresponding critical values in $[\eta, \beta)$ and
$$
m=\text{dim}(V\cap W)-\text{codim}(V\oplus W).
$$
\end{theorem}

\section{Proof of Theorem~\ref{thm05}}\label{Sec05}

Let the corresponding functional $J: H^1_0(\Omega)\to\mathbb{R}$ of \eqref{eq49} be defined by
\begin{equation}\label{eq66}
J(u)=\frac{a}{2}\|u\|^2+\frac{b}{4}\|u\|^4-\frac{\lambda}{2}\|u\|_2^2,\
u\in H^1_0(\Omega).
\end{equation}
Clearly, the $J\in C^1( H^1_0(\Omega), \bbr)$ and
$$
J'(u)v=(a+b\|u\|^2)\int_{\Omega}\nabla
u \nabla vdx -\lambda\int_{\Omega}uvdx,\ u, \ v\in H^1_0(\Omega).
$$
Therefore critical points of $J$ are  the weak
solutions of \eqref{eq49}.\vspace{6pt}

\medskip\par\noindent{\bf Proof of of Theorem~\ref{thm05}.}\quad

\noindent (I):\quad Suppose on the contrary that
Eq.~\eqref{eq49} has a nontrivial
solution $u\in H^1_0(\Omega)\setminus\{0\}$ for $\lambda\in (0,a\lambda_1]$, then we have
$$
0=a\|u\|^2+b\|u\|^4-\lambda\|u\|^2_2\geq
(a-\frac{\lambda}{\lambda_1})\|u\|^2+b\|u\|^4>0,
$$
this is impossible.

Before giving the proofs of (II) and (III), we shall prove the
following claim:

\noindent{\bf Claim:}\quad  If $\lambda\in (a\lambda_k,
a\lambda_{k+1}]$ with $k\geq 1$ then $u$ is a nontrivial solution
of \eqref{eq49} if and only if $u$ has the form of
\begin{equation}\label{eq65}
\overline{u}_j:=\sqrt{\frac{\lambda-a\lambda_j}{b\lambda_j^2}}\psi_j,
\end{equation}
where $\psi_j\in M(\lambda_j)$ satisfying
$$
\psi_j=\sum_{i=1}^{i_j}c_{ij}\varphi_{ij}\ \text{and}\
\sum_{i=1}^{i_j}c^2_{ij}=1, \quad j=1,2,\cdots k.
$$

Indeed, $u$ is a nontrivial solution of \eqref{eq49} if and only if
$u$ is a solution of the eigenvalue problem
$$
\left\{\aligned -\Delta u &=\frac{\lambda}{a+b\|u\|^2} u, &\quad \text{in}\ \Omega, \\
u&=0,& \text{on}\ \partial\Omega.
\endaligned
\right.
$$
Thus we must get
$$
\frac{\lambda}{a+b\|u\|^2}=\lambda_j \ \text{and}\ u=t_j\psi_j,
$$
for some $t_j>0$, $\psi_j\in M(\lambda_j)$ with $\|\psi_j\|_2^2=1$,
i.e., $\lambda=a\lambda_j+b\lambda_jt_j^2\|\psi_j\|^2\leq
a\lambda_{k+1}$, which implies that
$$
t_j=\sqrt{\frac{\lambda-a\lambda_j}{b\lambda_j^2}}, \ j=1,2,\cdots
k,
$$
so that
$$
u=\sqrt{\frac{\lambda-a\lambda_j}{b\lambda_j^2}}\psi_j, \ \psi_j\in
M(\lambda_j)\ \text{with}\
\psi_j=\sum_{i=1}^{i_j}c_{ij}\varphi_{ij}\ \text{and}\
\sum_{i=1}^{i_j}c^2_{ij}=1, \quad j=1,2,\cdots k.
$$
On the other hand, if $u$ has the form of \eqref{eq65}, then it is
easy to see that $u$ is a solution of \eqref{eq49}.  Therefore we
complete the proof of this claim. \vspace{6pt}

\noindent (II):\quad For the case of $\lambda\in (a\lambda_1, a\lambda_2]$, since  $\psi_1=\varphi_{11}$,
it follows from the above claim that \eqref{eq49} has the unique
positive solution
$$
\overline{u}_1=\sqrt{\frac{\lambda-a\lambda_1}{b\lambda_1^2}}\varphi_{11}.
$$

\noindent (III):\quad For the case $\lambda\in (a\lambda_k,
a\lambda_{k+1}]$, by the above claim, we see that \eqref{eq49}
has the solution of the form \eqref{eq65}.  We will complete the proof of Theorem~\ref{thm05} by going through the following two steps.

\noindent{\bf Step 1.} \quad $J(\overline{u}_1)=\min\{J(\overline{u}_j):
j=1,2,\cdots k\}$.

Indeed, for $j\in \{1,2,\cdots k\}$, we have
$$
J(\overline{u}_j)=\frac{a(\lambda-a\lambda_j)}{2b\lambda_j^2}\|\psi_j\|^2+
\frac{b}{4}\Big(\frac{\lambda-a\lambda_j}{b\lambda_j^2}\Big)^2\|\psi_j\|^4-
\frac{\lambda(\lambda-a\lambda_j)}{2b\lambda_j^2}\|\psi_j\|_2^2=-\frac{(\lambda-a\lambda_j)^2}{4b\lambda_j^2}.
$$
Set $g(t)=-(\lambda-at)^2/(4bt^2)$, $t\in(0, \lambda/a)$, by a
direct calculation, we see that $g'(t)>0$ for $t\in (0, \lambda/a)$,
therefore, $g(\lambda_1)<g(\lambda_2)<\cdots<g(\lambda_k)$
because $\lambda_1<\lambda_2<\cdots<\lambda_k$, so that
$J(\overline{u}_1)=\min\{J(\overline{u}_j): j=1,2,\cdots k\}$.

\noindent{\bf Step 2.} $\overline{u}_j$ is sign-changing,
$j=2,3,\cdots k$.

In fact, arguing by contradiction, without loss of generality, we may assume
that $\overline{u}_j\geq 0$ for some $2\leq j\leq k$.  Recalling
that $\varphi_{11}>0$ in $\Omega$ and $\overline{u}_j\in
M(\lambda_j)$ with
$$
\overline{u}_j=\sqrt{\frac{\lambda-a\lambda_j}{b\lambda_j^2}}\psi_j=\sqrt{\frac{\lambda-a\lambda_j}{b\lambda_j^2}}\sum_{i=1}^{i_j}c_{ij}\varphi_{ij}\
\text{and}\ \sum_{i=1}^{i_j}c^2_{ij}=1,
$$
then we obtain
$$
0\leq \int_{\Omega}\overline{u}_j\varphi_{11}dx=
\sqrt{\frac{\lambda-a\lambda_j}{b\lambda_j^2}}\sum_{i=1}^{i_j}c_{ij}\int_{\Omega}\varphi_{ij}\varphi_{11}dx=0,
$$
which means that $\overline{u}_j\equiv 0$, this is a contradiction since
$\|\overline{u}_j\|_2^2=(\lambda-a\lambda_j)/(b\lambda_j^2)>0$.
\qquad
\raisebox{-0.5mm}{\rule{1.5mm}{4mm}}\vspace{6pt}

\section{Proof of Theorems~\ref{thm01} and \ref{thm03}}\label{Sec03}

In this section, we are interested in the existence  of solutions to
\eqref{eq02} for $\lambda\geq a\lambda_1$.  We note that
the corresponding functional $I$ given by \eqref{eq04} with $q=2$ may be
indefinite.  We first establish a lemma which gives
the description of $(PS)_c$ sequence of $I$ for the above cases.
\begin{lemma}\label{lem01}
Suppose that $a$, $\lambda$, $\delta>0$, $b\geq0$ satisfy
$$
a\lambda_1\leq\lambda\leq
\frac{a\lambda_1\delta}{\delta-b\mathcal{S}^2}, \quad
\frac{b|\Omega|\lambda^2}{a^2}<\delta.
$$
If $\{u_n\}\subset H^1_0(\Omega)$ is a $(PS)_c$ sequence of $I$ with
$c\in(0, \frac{(a\mathcal{S})^2}{4(\delta-b\mathcal{S}^2)})$, then
$\{u_n\}$ strongly converges in $H^1_0(\Omega)$ up to a subsequence.
\end{lemma}
\begin{prooff}
Let $\{u_n\}$ be a $(PS)_c$ sequence for $I$ with $c\in(0,
\frac{(a\mathcal{S})^2}{4(\delta-b\mathcal{S}^2)})$, then
\begin{equation}\label{eq62}
I(u_n)\to c\ \text{and}\ I'(u_n)\to 0 \ \text{in}\ H^{-1}(\Omega)\
\text{as $n\to \infty$}.
\end{equation}
We claim that $\{u_n\}$ is bounded in $H^1_0(\Omega)$. If not,
arguing by contradiction, we may assume that  $\|u_n\|\to \infty$ as
$n\to \infty$.  Let $w_n=u_n/\|u_n\|$, we see that $\|w_n\|=1$, and
then there exists $w_0\in H^1_0(\Omega)$ such that as $n\to\infty$,
$$
w_n\rightharpoonup w_0 \ \text{in}\ H^1_0(\Omega), \ w_n\to w_0 \
\text{in}\ L^p(\Omega)\ (p\in(1, 2^*)).
$$
We divide into the following two cases to discussion.

{\bf Case $(i):$} $w_0=0$.

It follows from \eqref{eq62} and $\|u_n\|\to \infty$ that
$$
o(1)=\frac{I(u_n)-\frac14I'(u_n)u_n}{\|u_n\|^2}=
\frac{a}{4}-\frac{\lambda}{4}\|w_n\|_2^2 =\frac{a}{4}+o(1),
$$
this implies $a=0$, which contradicts $a>0$.

{\bf Case $(ii):$} $w_0\neq 0$, i.e., there exists $C_1>0$ such that
$\|w_n\|\geq \|w_0\|\geq C_1$, up to a subsequence.

Noting that $\|u_n\|\to\infty$ as $n\to\infty$, we see that
$$\aligned
o(1)=\frac{I(u_n)}{\|u_n\|^4}=&\frac{a}{2\|u_n\|^2}+\frac{b}{4}-\frac{\lambda\|u_n\|^2_2}{2\|u_n\|^4}
-\frac{\delta\|u_n\|^4_4}{4\|u_n\|^4} \\
=& \frac{b}{4}-\frac{\delta}{4}\|w_n\|_4^4+o(1),
\endaligned
$$
which means that
\begin{equation}\label{eq11}
b=\delta\|w_n\|^4_4+o(1).
\end{equation}
On the other hand, we also obtain that
$$
o(1)=\frac{I(u_n)-\frac14I'(u_n)u_n}{\|u_n\|^2}=
\frac{a}{4}-\frac{\lambda}{4}\|w_n\|_2^2
$$
this, together with \eqref{eq11}, yields that
$$
a+o(1)=\lambda\|w_n\|_2^2\leq
\lambda|\Omega|^{\frac12}\|w_n\|_4^2\leq
\lambda|\Omega|^{\frac12}\Big(\frac{b+o(1)}{\delta}\Big)^{\frac12},
$$
so that
$$
\delta\leq \frac{b\lambda^2|\Omega|}{a^2},
$$
this is a contradiction since $\delta>b\lambda^2|\Omega|/a^2$,
therefore our claim holds and $\{u_n\}$ is bounded.

We will complete the proof to show
that $\{u_n\}$ has a subsequence which strongly converges in
$H^1_0(\Omega)$. In fact, otherwise, by Proposition~\ref{pro01}, we
know that there exist a function $u_0\in
H^1_0(\Omega)$ which is a weak limit of $\{u_n\}$, a number $k\in
\bbn$ and further, for every $i\in\{1, 2,\cdots k\}$, a sequence of
values $\{R^i_n\}\subset \bbr^+$, points $\{x^i_n\}\subset
\overline{\Omega}$ and a  function $v_i\in
\mathcal{D}^{1,2}(\bbr^4)$ satisfying \eqref{eq07}, \eqref{eq08} and
\begin{equation}\label{eq13}
I(u_n)=\widetilde{I}(u_0)+\sum_{i=1}^k\widetilde{I}^{\infty}(v_i)+o(1)
\end{equation}
up to subsequences, where $o(1)\to 0$ as $n\to 0$ and
$\widetilde{I}$ and $\widetilde{I}^{\infty}$ are respectively given
by \eqref{eq09} and \eqref{eq10}.

It follows from the Poincar\'{e} inequality that
\begin{equation}\label{eq14}\aligned
\widetilde{I}(u_0)=&\widetilde{I}(u_0)-\frac{1}{4}\Big\{(a+bA)\|u_0\|^2-\lambda\|u_0\|_2^2-\delta\|u_0\|_4^4\Big\}\\
\geq &\frac{a}{4}(1-\frac{\lambda}{a\lambda_1})\|u_0\|^2,
\endaligned
\end{equation}
where $A:= \|u_0\|^2+\sum_{i=1}^k\|v_i\|^2_{\mathcal{D}}$.  By \eqref{eq10} and the Sobolev inequality we get
$$\aligned
0=&(a+bA)\|v_i\|^2_{\mathcal{D}}-\delta\int_{\bbr^4}|v_i|^4dx \\
\geq &
\big(a+b\|u_0\|^2+b\|v_i\|^2_{\mathcal{D}}\big)\|v_i\|^2_{\mathcal{D}}-\frac{\delta\|v_i\|^4_{\mathcal{D}}}{\mathcal{S}^2}
\\
= &
(a+b\|u_0\|^2)\|v_i\|^2_{\mathcal{D}}-\frac{\delta-b\mathcal{S}^2}{\mathcal{S}^2}\|v_i\|^4_{\mathcal{D}},
\endaligned
$$
so that
\begin{equation}\label{eq15}
\|v_i\|^2_{\mathcal{D}}\geq
\frac{\mathcal{S}^2(a+b\|u_0\|^2)}{\delta-b\mathcal{S}^2},
\end{equation}
here we have used the inequality $b\mathcal{S}^2<\delta$ which can be obtained from the
assumptions $\lambda\geq a\lambda_1$,
$\delta>b|\Omega|\lambda^2/a^2$ and the following inequality
\begin{equation}\label{eq17}
0<\|\varphi_{11}\|^4=\lambda^2_1\|\varphi_{11}\|_2^4\leq
\frac{\lambda^2}{a^2}\|\varphi_{11}\|_2^4\leq
\frac{\lambda^2|\Omega|}{a^2}\|\varphi_{11}\|_4^4<
\frac{\lambda^2|\Omega|}{a^2\mathcal{S}^2}\|\varphi_{11}\|^4,
\end{equation}
where $\varphi_{11}$ is the positive eigenfunction of the first
eigenvalue $\lambda_1$ of $(-\Delta, H_0^1(\Omega))$.

From \eqref{eq08} and \eqref{eq10}, we see that for
$i\in\{1,2,\cdots,k\}$,
\begin{equation}\label{eq16}
\widetilde{I}^{\infty}(v_i)=\widetilde{I}^{\infty}(v_i)-\frac{1}{4}\Big\{(a+bA)\|v_i\|^2_{\mathcal{D}}-\int_{\bbr^4}|v_i|^4dx\Big\}\geq
\frac{a}{4}\|v_i\|^2_{\mathcal{D}},
\end{equation}
which, together with \eqref{eq13}, \eqref{eq14}, \eqref{eq15} and
\eqref{eq16}, implies that
$$
\aligned
c=\lim_{n\to\infty}I(u_n)=&\widetilde{I}(u_0)+\sum_{i=1}^k\widetilde{I}^{\infty}(v_i)\\
\geq &
\frac{a}{4}(1-\frac{\lambda}{a\lambda_1})\|u_0\|^2+\frac{(a\mathcal{S})^2}{4(\delta-b\mathcal{S}^2)}
+ \frac{ab\mathcal{S}^2\|u_0\|^2}{4(\delta-b\mathcal{S}^2)} \\
=& \frac{(a\mathcal{S})^2}{4(\delta-b\mathcal{S}^2)} +
\frac{a}{4}\Big(\frac{b\mathcal{S}^2}{\delta-b\mathcal{S}^2}+1-\frac{\lambda}{a\lambda_1}\Big)\|u_0\|^{2}\\
=&\frac{(a\mathcal{S})^2}{4(\delta-b\mathcal{S}^2)} + \frac{a}{4}\Big(\frac{\delta}{\delta-b\mathcal{S}^2}-\frac{\lambda}{a\lambda_1}\Big)\|u_0\|^{2}\\
\geq & \frac{(a\mathcal{S})^2}{4(\delta-b\mathcal{S}^2)}
\endaligned
$$
since $\lambda\leq a\lambda_1\delta/(\delta-b\mathcal{S}^2)$.  This
contradicts $0<c<(a\mathcal{S})^2/(4(\delta-b\mathcal{S}^2))$, and then we complete proof of this lemma.
\end{prooff}
\begin{remark}\label{rmk04}  {\em Our case of $\lambda\geq
a\lambda_1$ is quite different from the case of
$\lambda\in(0, a\lambda_1)$ assumed in \cite{D14} since, except for the indefiniteness of the operator $-\Delta-a\lambda$,  it is not obvious whether
$\widetilde{I}(u_0)\geq 0$ or not from \eqref{eq14}, and now it follows from
\eqref{eq17} that
\begin{equation}\label{eq3.9}
\mathcal{S}^2<\frac{\lambda^2|\Omega|}{a^2},
\end{equation}
which is also different from the
assumption given by \cite[Lemma 2.1]{D14}.}
\end{remark}

To prove Theorems~\ref{thm01} and \ref{thm03}, we will find two suitable closed
subspace $V$ and $W$ with $V\cap W\neq \{0\}$ and $V\oplus W=H$ such
that the functional $I$ satisfies the condition $(a_1)-(a_3)$ of
Theorem~\ref{thm02} with
$\beta=(a\mathcal{S})^2/(4(\delta-b\mathcal{S}^2))$.

Given $a$, $\lambda$, $\rho>0$, we set
$$
\lambda^+=\min\{\lambda_j: \lambda<a\lambda_j\}, \quad
\lambda^-=\max\{\lambda_j: a\lambda_j\leq \lambda\},
$$
\begin{equation}\label{eq19}
H_1=\overline{\oplus_{\lambda_j\geq \lambda^+}M(\lambda_j)}, \quad
H_2=\oplus_{\lambda_j\leq \lambda^-}M(\lambda_j),\quad B_{\rho}=\{x\in \bbr^4: |x|<\rho\}.
\end{equation}

Without loss of generality, we can suppose that $0\in \Omega$ and
that $B_{1}(0)\subset \Omega$. Given $\varepsilon>0$, let
$U_{\varepsilon}(x)$ be the function given by \eqref{eq20}, we
define
$$
\Psi_{\varepsilon}(x)=\phi(x)U_{\varepsilon}(x),
$$
where $\phi\in C^{\infty}_0(B_1(0))$ with $\phi(x)=1$ in
$B_{1/2}(0)$.  Similarly as in \cite{BN83,CFP85}, we have the following energy estimates of $\Psi_{\varepsilon}$:
\begin{lemma}\label{lem02}
\begin{equation}\label{eq22}
\|\Psi_{\varepsilon}\|^2=\mathcal{S}^2+ O(\varepsilon);
\end{equation}
\begin{equation}\label{eq23}
\|\Psi_{\varepsilon}\|^4_4=\mathcal{S}^2+ O(\varepsilon^2);
\end{equation}
\begin{equation}\label{eq24}
\|\Psi_{\varepsilon}\|_2^2\geq K_1\varepsilon|\log\varepsilon|+
O(\varepsilon),
\end{equation}
\begin{equation}\label{eq25}
\|\Psi_{\varepsilon}\|_1\leq K_2\varepsilon^{\frac12},
\end{equation}
\begin{equation}\label{eq26}
\|\Psi_{\varepsilon}\|_3^3\leq K_3\varepsilon^\frac{1}{2}
\end{equation}
\end{lemma}

Set $\overline{j}=\max\{j: \lambda_j\leq\lambda\}$ and denote by
$P_j$ the orthogonal projector onto the eigenspace $M(\lambda_j)$
corresponding to $\lambda_j$. We now set
\begin{equation}\label{eq27}
\widetilde{\Psi}_{\varepsilon}=\Psi_{\varepsilon}-\sum_{j=1}^{\overline{j}}P_j\Psi_{\varepsilon}
\end{equation}
and give some energy estimates for $\widetilde{\Psi}_{\varepsilon}$ in the following lemma.

\begin{lemma}\label{lem03}
For $\varepsilon>0$ sufficiently small, we have
\begin{equation}\label{eq28}
\|\sum_{j=1}^{\overline{j}}P_j\Psi_{\varepsilon}\|_2^2\leq
C\varepsilon;
\end{equation}
\begin{equation}\label{eq31}
\|\sum_{j=1}^{\overline{j}}P_j\Psi_{\varepsilon}\|_{\infty}\leq
C\varepsilon^{\frac12};
\end{equation}
\begin{equation}\label{eq29}
\frac{a\|\widetilde{\Psi}_{\varepsilon}\|^2-\lambda\|\widetilde{\Psi}_{\varepsilon}\|_2^2}{\|\widetilde{\Psi}_{\varepsilon}\|_4^2}=
a\mathcal{S}+C\varepsilon\log\varepsilon+ O(\varepsilon),
\end{equation}
where $C>0$ denote variant different constants.
\end{lemma}
\begin{prooff}
Since $\{\varphi_{ij}\}$ be the orthogonal eigenfunctions
corresponding to $\lambda_j$ spanning $M(\lambda_j)$, $i=1,2,\cdots
i_j$, where $i_j=\text{dim}M(\lambda_j)$, then
\begin{equation}\label{eq30}
\int_{\Omega}\varphi_{ij}\varphi_{kn}dx=\left\{\aligned &1, \quad
\text{if}\ i=k\ \text{and}\ j=n, \\
&0, \quad \text{others}.
\endaligned
\right.
\end{equation}
Since $\Psi_{\varepsilon}\in H^1_0(\Omega)$ and
$\{\varphi_{ij}\}$($i=1,2,\cdots i_j$; $j=1, 2,\cdots$) is a family
of orthonormal basis of $H^1_0(\Omega)$, there holds
\begin{equation}\label{eq33}
\Psi_{\varepsilon}(x)=\sum_{j=1}^{\infty}\sum_{i=1}^{i_j}\Big(\int_{\Omega}\Psi_{\varepsilon}\varphi_{ij}dx\Big)\varphi_{ij}(x),
\end{equation}
then we have
\begin{equation}\label{eq34}
P_j\Psi_{\varepsilon}(x)=\sum_{i=1}^{i_j}\Big(\int_{\Omega}\Psi_{\varepsilon}\varphi_{ij}dx\Big)\varphi_{ij}(x).
\end{equation}
On the other hand, recalling that $\varphi_{ij}\in
C^{\infty}(\Omega)$, we see that there exists $C>0$, dependent of
$j$, such that $\sup\{\varphi_{ij}(x): x\in \Omega, i=1,2,\cdots
i_j\}\leq C$, therefore, by \eqref{eq25} and \eqref{eq34} we know
\begin{equation}\label{eq37}
\|P_j\Psi_{\varepsilon}\|_2^2=\sum_{i=1}^{i_j}\Big(\int_{\Omega}\Psi_{\varepsilon}\varphi_{ij}dx\Big)^2\leq
C\big(\sum_{i=1}^{i_j}\|\varphi_{ij}\|^2_{\infty}\big)\|\Psi_{\varepsilon}\|^2_1\leq
C\|\Psi_{\varepsilon}\|^2_1\leq C\varepsilon.
\end{equation}
Similarly, we also deduce
\begin{equation}\label{eq38}
\|\sum_{j=1}^{\overline{j}}P_j\Psi_{\varepsilon}\|_2^2=\sum_{j=1}^{\overline{j}}\|P_j\Psi_{\varepsilon}\|_2^2\leq
C\varepsilon,
\end{equation}
$$
\|P_j\Psi_{\varepsilon}\|_{\infty}\leq
\sum_{i=1}^{i_j}\Big(\int_{\Omega}\Psi_{\varepsilon}\varphi_{ij}dx\Big)\|\varphi_{ij}(x)\|_{\infty}\leq
C\|\Psi_{\varepsilon}\|_1\leq C\varepsilon^{\frac12},
$$
$$
\|\sum_{j=1}^{\overline{j}}P_j\Psi_{\varepsilon}\|_{\infty}=\sum_{j=1}^{\overline{j}}\|P_j\Psi_{\varepsilon}\|_{\infty}\leq
C\varepsilon^{\frac12}.
$$
Thus \eqref{eq28} and \eqref{eq31} hold.  Furthermore, we have, by
\eqref{eq26} and \eqref{eq38},
$$\aligned
\Big|\int_{\Omega}(|\widetilde{\Psi}_{\varepsilon}|^4-|\Psi_{\varepsilon}|^4)dx\Big|=&
4\Big|\int_{\Omega}\int_0^1\Big|\Psi_{\varepsilon}-t\sum_{j=1}^{\overline{j}}P_j\Psi_{\varepsilon}\Big|^{2}\Big(\Psi_{\varepsilon}
-t\sum_{j=1}^{\overline{j}}P_j\Psi_{\varepsilon}\Big)\Big(\sum_{j=1}^{\overline{j}}P_j\Psi_{\varepsilon}\Big)dtdx\Big|
\\
\leq &32\Big|\int_{\Omega}\int_0^1\Big(|\Psi_{\varepsilon}|^{3}+
\Big|t\sum_{j=1}^{\overline{j}}P_j\Psi_{\varepsilon}\Big|^3\Big)\Big(\sum_{j=1}^{\overline{j}}P_j\Psi_{\varepsilon}\Big)dtdx\Big|
\\
\leq &
32\Big(\|\Psi_{\varepsilon}\|^3_3\|\sum_{j=1}^{\overline{j}}P_j\Psi_{\varepsilon}\|_{\infty}+\|\sum_{j=1}^{\overline{j}}P_j\Psi_{\varepsilon}\|^4_4\Big)\\
\leq & 32\Big(\|\Psi_{\varepsilon}\|^3_3\|\sum_{j=1}^{\overline{j}}P_j\Psi_{\varepsilon}\|_{\infty}+
\|\sum_{j=1}^{\overline{j}}P_j\Psi_{\varepsilon}\|^2_2\|\sum_{j=1}^{\overline{j}}P_j\Psi_{\varepsilon}\|^2_{\infty}\Big)\\
\leq & C(\varepsilon+\varepsilon^2),
\endaligned
$$
which implies that
\begin{equation}\label{eq32}
\|\widetilde{\Psi}_{\varepsilon}\|_4^4=\|\Psi_{\varepsilon}\|_4^4+O(\varepsilon).
\end{equation}
It follows from \eqref{eq30}, \eqref{eq33}, \eqref{eq34},
\eqref{eq37} and \eqref{eq38} that
\begin{equation}\label{eq35}\aligned
\|\widetilde{\Psi}_{\varepsilon}\|^2=&
\int_{\Omega}\Big|\nabla\Psi_{\varepsilon}-\sum_{j=1}^{\overline{j}}\nabla\big(P_j\Psi_{\varepsilon}\big)\Big|^2dx
\\
=&\int_{\Omega}\Big\{|\nabla\Psi_{\varepsilon}|^2-2\sum_{j=1}^{\overline{j}}\nabla\Psi_{\varepsilon}\nabla\big(P_j\Psi_{\varepsilon}\big)+
\sum_{j=1}^{\overline{j}}\Big|\nabla\big(P_j\Psi_{\varepsilon}\big)\Big|^2\Big\}dx
\\
=&\|\Psi_{\varepsilon}\|^2-2\sum_{j=1}^{\overline{j}}\sum_{i=1}^{i_j}\int_{\Omega}\Psi_{\varepsilon}\varphi_{ij}dx\int_{\Omega}\nabla\Psi_{\varepsilon}\nabla\varphi_{ij}dx
+\sum_{j=1}^{\overline{j}}\sum_{i=1}^{i_j}\Big(\int_{\Omega}\Psi_{\varepsilon}\varphi_{ij}dx\Big)^2\int_{\Omega}|\nabla\varphi_{ij}|^2dx\\
=&\|\Psi_{\varepsilon}\|^2-\sum_{j=1}^{\overline{j}}\sum_{i=1}^{i_j}\lambda_j\Big(\int_{\Omega}\Psi_{\varepsilon}\varphi_{ij}dx\Big)^2\\
=&\|\Psi_{\varepsilon}\|^2-\sum_{j=1}^{\overline{j}}\lambda_j\|P_j\Psi_{\varepsilon}\|_2^2,
\endaligned
\end{equation}
\begin{equation}\label{eq36}\aligned
\|\widetilde{\Psi}_{\varepsilon}\|^2_2=&\int_{\Omega}\Big|\Psi_{\varepsilon}-\sum_{j=1}^{\overline{j}}\big(P_j\Psi_{\varepsilon}\big)\Big|^2dx\\
=&\int_{\Omega}\Big\{|\Psi_{\varepsilon}|^2-2\sum_{j=1}^{\overline{j}}\Psi_{\varepsilon}\big(P_j\Psi_{\varepsilon}\big)+
\sum_{j=1}^{\overline{j}}\Big|\big(P_j\Psi_{\varepsilon}\big)\Big|^2\Big\}dx\\
&\|\Psi_{\varepsilon}\|^2_2-2\sum_{j=1}^{\overline{j}}\sum_{i=1}^{i_j}\Big(\int_{\Omega}\Psi_{\varepsilon}\varphi_{ij}dx\Big)^2
+\sum_{j=1}^{\overline{j}}\sum_{i=1}^{i_j}\Big(\int_{\Omega}\Psi_{\varepsilon}\varphi_{ij}dx\Big)^2\int_{\Omega}|\varphi_{ij}|^2dx\\
=&\|\Psi_{\varepsilon}\|^2_2-\sum_{j=1}^{\overline{j}}\sum_{i=1}^{i_j}\Big(\int_{\Omega}\Psi_{\varepsilon}\varphi_{ij}dx\Big)^2\\
=&\|\Psi_{\varepsilon}\|^2_2-\sum_{j=1}^{\overline{j}}\|P_j\Psi_{\varepsilon}\|_2^2.
\endaligned
\end{equation}
Combining \eqref{eq32},\eqref{eq35},\eqref{eq36} and
\eqref{eq22}-\eqref{eq24} of Lemma~\ref{lem02}, we obtain that
$$\aligned
\frac{a\|\widetilde{\Psi}_{\varepsilon}\|^2-\lambda\|\widetilde{\Psi}_{\varepsilon}\|_2^2}{\|\widetilde{\Psi}_{\varepsilon}\|_4^2}=&
\frac{a\|\Psi_{\varepsilon}\|^2-\lambda\|\Psi_{\varepsilon}\|_2^2+\sum_{j=1}^{\overline{j}}(\lambda-a\lambda_j)\|P_j\Psi_{\varepsilon}\|_2^2}
{\|\Psi_{\varepsilon}\|_4^2+O(\varepsilon)} \\
=&
\frac{a\|\Psi_{\varepsilon}\|^2-\lambda\|\Psi_{\varepsilon}\|_2^2+O(\varepsilon)}
{\|\Psi_{\varepsilon}\|_4^2+O(\varepsilon^2)} \\
=&\frac{a\|\Psi_{\varepsilon}\|^2-\lambda\|\Psi_{\varepsilon}\|_2^2}
{\|\Psi_{\varepsilon}\|_4^2}+O(\varepsilon)\\
=&a\mathcal{S}+C\varepsilon\log\varepsilon+ O(\varepsilon),
\endaligned
$$
then \eqref{eq29} is proved and this completes the proof of the
lemma.
\end{prooff}

Let $H_1$ and $H_2$ be given in \eqref{eq19}.  Clearly $\widetilde{\Psi}_{\varepsilon}$ defined by \eqref{eq27} belongs to $H_1$.  Now we define a subspace $\overline{W}_{\varepsilon}$ in
$H^1_0(\Omega)$ as
\begin{equation}\label{eq39}
\overline{W}_{\varepsilon}=\{u\in H^1_0(\Omega):
u=\overline{u}+t\widetilde{\Psi}_{\varepsilon}, \overline{u}\in H_2,
t\in \bbr\},
\end{equation}
then we have the following result:
\begin{lemma}\label{lem04}
Suppose that $a$, $\lambda>0$, then for $\varepsilon>0$ sufficiently
small,
$$
\sup_{u\in
\overline{W}_{\varepsilon}\setminus\{0\}}\frac{\|u\|^4}{\|u\|_4^4}\leq
2|\Omega|\Big(\sum_{j=1}^{\overline{j}}i_j\Big)^2\frac{\lambda^2}{a^2}.
$$
\end{lemma}
\begin{prooff}
By the definition of $\overline{W}_{\varepsilon}$ and the orthogonality of $H_1$ and $H_2$, we know that for each $u\in \overline{W}_{\varepsilon}$,
$u=\overline{u}+t\widetilde{\Psi}_{\varepsilon}$ for some $\overline{u}\in H_2$ and
$t\in \bbr$, and
\begin{equation}\label{eq40}
\|u\|^4=\|\overline{u}+t\widetilde{\Psi}_{\varepsilon}\|^4=(\|\overline{u}+t\widetilde{\Psi}_{\varepsilon}\|^2)^2\leq
2\|\overline{u}\|^4+2\|t\widetilde{\Psi}_{\varepsilon}\|^4.
\end{equation}
Since
$$
\|\overline{u}\|_{\infty}\leq C, \quad
\|\widetilde{\Psi}_{\varepsilon}\|_{\infty}=\|\Psi_{\varepsilon}-\sum_{j=1}^{\overline{j}}P_j\Psi_{\varepsilon}\|_{\infty}\leq
\|\Psi_{\varepsilon}\|_{\infty}+\|\sum_{j=1}^{\overline{j}}P_j\Psi_{\varepsilon}\|_{\infty}\leq
8\varepsilon^{-\frac{1}{2}}+C\varepsilon^{\frac{1}{2}},
$$
we have
$$
\Big|\int_{\Omega}\overline{u}(t\widetilde{\Psi}_{\varepsilon})^3dx\Big|\leq
|t|^3\|\widetilde{\Psi}_{\varepsilon}\|_{\infty}\Big|\int_{\Omega}\overline{u}\widetilde{\Psi}_{\varepsilon}dx\Big|=0,
\quad
\int_{\Omega}\overline{u}^3(t\widetilde{\Psi}_{\varepsilon})dx=\int_{\Omega}\overline{u}^2(t\widetilde{\Psi}_{\varepsilon})^2dx=0,
$$
so that
\begin{equation}\label{eq41}\aligned
\|u\|^4_4=&\int_{\Omega}|\overline{u}+t\widetilde{\Psi}_{\varepsilon}|^4dx\\
=&\int_{\Omega}(|\overline{u}|^4+4u^3(t\widetilde{\Psi}_{\varepsilon})
+6u^2(t\widetilde{\Psi}_{\varepsilon})^2+4u(t\widetilde{\Psi}_{\varepsilon})^3+|t\widetilde{\Psi}_{\varepsilon}|^4)dx\\
=&\|\overline{u}\|_4^4+\|t\widetilde{\Psi}_{\varepsilon}\|_4^4.
\endaligned
\end{equation}
Therefore, by using \eqref{eq40}, \eqref{eq41} and
Lemma~\ref{lem05}, we see that
\begin{equation}\label{eq42}\aligned
\sup_{u\in
\overline{W}_{\varepsilon}\setminus\{0\}}\frac{\|u\|^4}{\|u\|_4^4}\leq
&\sup_{u\in\overline{W}_{\varepsilon}\setminus\{0\}}\frac{2\|\overline{u}\|^4+2\|t\widetilde{\Psi}_{\varepsilon}\|^4}
{\|\overline{u}\|_4^4+\|t\widetilde{\Psi}_{\varepsilon}\|^4_4}
\\ \leq& 2\max\Big\{\sup_{\overline{u}\in H_2\setminus\{0\}}\frac{\|\overline{u}\|^4}{\|\overline{u}\|_4^4},
\quad \frac{\|t\widetilde{\Psi}_{\varepsilon}\|^4}{\|t\widetilde{\Psi}_{\varepsilon}\|^4_4}\Big\} \\
=& 2\max\Big\{\sup_{\overline{u}\in
H_2\setminus\{0\}}\frac{\|\overline{u}\|^4}{\|\overline{u}\|_4^4},
\quad
\frac{\mathcal{S}^4+O(\varepsilon)}{\mathcal{S}^2+O(\varepsilon^2)}\Big\}.
\endaligned
\end{equation}
On the other hand, for each $\overline{u}\in H_2\setminus\{0\}$,
there exists $\{c_{ij}\}\subset\mathbb{R}\ (j=1,2,\cdots,
\overline{j}; \ i=1,2,\cdots, i_j)$ such that
$$
\overline{u}=\sum_{j=1}^{\overline{j}}\sum_{i=1}^{i_j}c_{ij}\varphi_{ij},
$$
where $\{\varphi_{ij}\}$ is the family of orthonormal basis of
$H^1_0(\Omega)$ given in the proof of Lemma~\ref{lem03}. Therefore
we have
$$
\|\overline{u}\|^4=\Big(\sum_{j=1}^{\overline{j}}\sum_{i=1}^{i_j}\|c_{ij}\varphi_{ij}\|^2\Big)^2\leq
\Big(\sum_{j=1}^{\overline{j}}i_j\Big)^2\Big(\sum_{j=1}^{\overline{j}}\sum_{i=1}^{i_j}\|c_{ij}\varphi_{ij}\|^4\Big),
$$
$$
\|\overline{u}\|^4_4=\sum_{j=1}^{\overline{j}}\sum_{i=1}^{i_j}\|c_{ij}\varphi_{ij}\|^4_4.
$$
By using Lemma~\ref{lem05} again, we get
\begin{equation}\label{eq43}
\sup_{\overline{u}\in
H_2\setminus\{0\}}\frac{\|\overline{u}\|^4}{\|\overline{u}\|_4^4}\leq
\Big(\sum_{j=1}^{\overline{j}}i_j\Big)^2
\max\Big\{\frac{\|\varphi_{ij}\|^4}{\|\varphi_{ij}\|^4_4},
j=1,2,\cdots \overline{j}, \quad i=1,2,\cdots, i_j\Big\}.
\end{equation}
It follows from \eqref{eq17}, \eqref{eq42}, \eqref{eq43} and \eqref{eq3.9} in
Remark~\ref{rmk04} that
$$\aligned
\sup_{u\in
\overline{W}_{\varepsilon}\setminus\{0\}}\frac{\|u\|^4}{\|u\|_4^4}\leq&
2\max\Big\{\frac{\mathcal{S}^4+O(\varepsilon)}{\mathcal{S}^2+O(\varepsilon^2)},
\quad \Big(\sum_{j=1}^{\overline{j}}i_j\Big)^2
\frac{\|\varphi_{ij}\|^4}{\|\varphi_{ij}\|^4_4},\quad j=1,2,\cdots
\overline{j}, \quad i=1,2,\cdots, i_j\Big\}\\
\leq&
2\max\Big\{\frac{\mathcal{S}^4+O(\varepsilon)}{\mathcal{S}^2+O(\varepsilon^2)},
\quad \Big(\sum_{j=1}^{\overline{j}}i_j\Big)^2(\lambda_j)^2|\Omega|,
\quad j=1,2,\cdots \overline{j}\Big\}\\
\leq &
2|\Omega|\Big(\sum_{j=1}^{\overline{j}}i_j\Big)^2\frac{\lambda^2}{a^2}
\endaligned
$$
for $\varepsilon>0$ sufficiently small.
\end{prooff}

\begin{lemma}\label{lem06}
Suppose that $a>0$ and $\lambda\in[a\lambda_1,
a\lambda_1\delta/(\delta-b\mathcal{S}^2))$, then for $\varepsilon>0$
small,
\begin{equation}\label{eq63}
\sup_{u\in
\overline{W}_{\varepsilon}}I(u)<\frac{(a\mathcal{S})^2}{4(\delta-b\mathcal{S}^2)}.
\end{equation}
provided one of the following conditions holds
\begin{itemize}
  \item [(1)] $\delta>0$, $b>0$ is sufficiently small;
  \item [(2)] $b>0$, $\delta>0$ is sufficiently large.
\end{itemize}
\end{lemma}
\begin{prooff}
First of all, let $\delta>0$ be fixed, we may choose $b_1>0$ such that for each $b\in (0,
b_1)$, it holds
\begin{equation}\label{eq3.36}
\frac{a\lambda_1\delta}{\delta-b\mathcal{S}^2}\leq a\lambda_2, \quad
4b\lambda_2|\Omega|<\delta,
\end{equation}
thus $H_2=M(\lambda_1)$.  It follows from Lemma~\ref{lem04} that there exists
$\varepsilon_0>0$ such that for $\varepsilon\in (0, \varepsilon_0)$
$$
\sup_{u\in\overline{W}_{\varepsilon}\setminus\{0\}}\frac{\|u\|^4}{\|u\|_4^4}\leq
2\lambda^2_2|\Omega|.
$$
For each $u\in\overline{W}_{\varepsilon}$, we we may assume that
$a\|u\|^2-\lambda\|u\|^2>0$ (since otherwise $I(tu)\leq 0$ for all
$t>0$, which implies $\sup_{u\in \overline{W}_{\varepsilon}}I(u)\leq
\sup_{u\in \overline{W}_{\varepsilon}}\max_{t>0}I(tu)\leq 0$, so
that the desired inequality \eqref{eq63} holds trivially), hence
$I(tu)>0$ for $t>0$ small, and
\begin{equation}\label{eq45}
0<\max_{t>0}I(tu)\leq
\frac{(a\|u\|^2-\lambda\|u\|_2^2)^2}{4(\delta-b\frac{\|u\|^4}{\|u\|_4^4})\|u\|_4^4}
\leq
\frac{(a\|u\|^2-\lambda\|u\|_2^2)^2}{4(\delta-2b\lambda_2^2|\Omega|)\|u\|_4^4}.
\end{equation}
Next, we claim that
$$
\sup\{a\|u\|^2-\lambda\|u\|_2^2:\ u\in\overline{W}_{\varepsilon},
\|u\|_4=1\}\leq
a\mathcal{S}+C\varepsilon\log\varepsilon+O(\varepsilon)
$$
holds for some positive constant $C$  and small $\varepsilon$.

Indeed, for $u\in\overline{W}_\varepsilon$ with $\|u\|_4=1$, we have
$$
1=\|u\|_4^4=\|\overline{u}\|_4^4+\|t\widetilde{\Psi}_{\varepsilon}\|_4^4,
$$
then by \eqref{eq32}, we see that $t$ is bounded and
$\|t\widetilde{\Psi}_{\varepsilon}\|_4^4\leq 1$. Therefore, by using
Lemma~\ref{lem03},
\begin{equation}\label{eq46}\aligned
a\|u\|^2-\lambda\|u\|_2^2=&a(\|\overline{u}\|^2+\|t\widetilde{\Psi}_{\varepsilon}\|^2_2)-\lambda(\|\overline{u}\|^2_2+\|t\widetilde{\Psi}_{\varepsilon}\|_2^2)\\
= &
a\|\overline{u}\|^2-\lambda\|\overline{u}\|^2_2+a\|t\widetilde{\Psi}_{\varepsilon}\|^2-\lambda\|t\widetilde{\Psi}_{\varepsilon}\|_2^2
\\ = & (a\lambda_1-\lambda)\|\overline{u}\|^2_2 +\frac{a\|\widetilde{\Psi}_{\varepsilon}\|^2-
\lambda\|\widetilde{\Psi}_{\varepsilon}\|_2^2}{\|\widetilde{\Psi}_{\varepsilon}\|_4^2}\|t\widetilde{\Psi}_{\varepsilon}\|_4^2\\
\leq& \frac{a\|\widetilde{\Psi}_{\varepsilon}\|^2-
\lambda\|\widetilde{\Psi}_{\varepsilon}\|_2^2}{\|\widetilde{\Psi}_{\varepsilon}\|_4^2}
\\ \leq & a\mathcal{S}+C\varepsilon\log\varepsilon+O(\varepsilon).
\endaligned
\end{equation}
It follows from \eqref{eq45} and \eqref{eq46} that
$$\aligned
\sup_{u\in \overline{W}_{\varepsilon}}I(u)\leq &
\frac{(a\mathcal{S})^2+C\varepsilon\log\varepsilon+O(\varepsilon)}{4(\delta-2b\lambda_2^2|\Omega|)}
\\ =& \frac{a^2}{4}\Big(\frac{\mathcal{S}^2}{\delta-b\mathcal{S}^2}+
\frac{b\mathcal{S}^2(2\lambda_2^2|\Omega|-\mathcal{S}^2)}{(\delta-b\mathcal{S}^2)(\delta-2b\lambda_2^2|\Omega|)}
+\frac{C\varepsilon\log\varepsilon+O(\varepsilon)}{\delta-2b\lambda_2^2|\Omega|}\Big)
\\ \leq & \frac{a^2}{4}\Big(\frac{\mathcal{S}^2}{\delta-b\mathcal{S}^2}+
\frac{4b\mathcal{S}^2(2\lambda_2^2|\Omega|-\mathcal{S}^2)}{\delta^2}
+\frac{C\varepsilon\log\varepsilon+O(\varepsilon)}{\delta}\Big),
\endaligned
$$
so that there exists $\varepsilon_1>0$ such that for each
$\varepsilon\in (0, \varepsilon_1)$,  we
can find $b_0:=b_0(\varepsilon)$ satisfying that for each $b\in(0, b_0)$ the inequality \eqref{eq63} remains true.

Finally, if $b>0$ is fixed then  \eqref{eq3.36} obviously holds for all large $\delta>0$.  Similarly as in the above, we see that
\eqref{eq63} also holds for all large $\delta>0$.  Thus we have proved this lemma.
\end{prooff}
\begin{remark}\label{rmk03}{\em
In Lemma~\ref{lem06}, we have shown that the assumption $(a_{31})$ of
Theorem~\ref{thm02} holds with $W=\overline{W}_{\varepsilon}$ and
$\beta=(a\mathcal{S})^2/(4(\delta-b\mathcal{S}^2))$, the choice of
$\beta$ is consistent with Lemma~\ref{lem01}.  However, by checking
the proof of Lemma~\ref{lem06} carefully, we observe that the
conclusion of Lemma~\ref{lem06} is still true if we take the place of
$(a\mathcal{S})^2/(4(\delta-b\mathcal{S}^2))$ by
$(a\mathcal{S})^2/(4\delta)-\sigma$ for some $\sigma>0$ small.
Indeed, it follows from \eqref{eq45} and \eqref{eq46} that
$$\aligned
\sup_{u\in \overline{W}_{\varepsilon}}I(u)\leq &
\frac{(a\mathcal{S})^2+C\varepsilon\log\varepsilon+O(\varepsilon)}{4(\delta-2b\lambda_2^2|\Omega|)}
\\ =& \frac{a^2}{4}\Big(\frac{\mathcal{S}^2}{\delta}+
\frac{2b\mathcal{S}^2\lambda_2^2|\Omega|}{\delta(\delta-2b\lambda_2^2|\Omega|)}
+\frac{C\varepsilon\log\varepsilon+O(\varepsilon)}{\delta-2b\lambda_2^2|\Omega|}\Big)
\\ \leq & \frac{a^2}{4}\Big(\frac{\mathcal{S}^2}{\delta}+
\frac{4b\mathcal{S}^2\lambda_2^2|\Omega|}{\delta^2}
+\frac{C\varepsilon\log\varepsilon+O(\varepsilon)}{\delta}\Big)\\
<& \frac{(a\mathcal{S})^2}{4\delta}-\sigma
\endaligned
$$
since for $\varepsilon$ sufficiently small, we can take
$b\ll\varepsilon$ or $\delta$ large enough.}
\end{remark}

\begin{lemma}\label{lem10}
Suppose that $a$, $\delta>0$.  Let $\{b_n\}$ and
$\{\overline{\lambda}_n\}$ be two sequences of positive numbers
satisfying $b_n\to 0$ and $\overline{\lambda}_n\to a\lambda_1$ as
$n\to \infty$, where $\overline{\lambda}_n>a\lambda_1$ for all $n\in
\bbn$.  If $u_n$ is the solution of Eq.~\eqref{eq02} corresponding
to $b_n$ and $\overline{\lambda}_n$ and $\{I(u_n)\}$ is bounded,
then $\{u_n\}$ is bounded.
\end{lemma}
\begin{prooff}
Arguing by contradiction, we may assume $\|u_n\|\to\infty$ as $n\to\infty$.  Let $w_n=u_n/\|u_n\|$,
it is easy to see that $\|w_n\|=1$, and then there exists $w_0\in
H^1_0(\Omega)$ such that $w_n\rightharpoonup w_0$ in $H^1_0(\Omega)$
and $w_n\to w_0$ in $L^p(\Omega)$ $(1<p<4)$ as $n\to\infty$, thus we
have
$$
0=\frac{I'(u_n)u_n}{\|u_n\|^4}=\frac{a}{\|u_n\|^2}+b_n-\frac{\overline{\lambda}_n\|w_n\|^2_2}{\|u_n\|^2}-\delta\|w_n\|^4_4,
$$
which implies that $\|w_n\|_4^4\to 0$ as $n\to\infty$, therefore we
obtain $w_0=0$.  On the other hand,
$$
o(1)=\frac{I(u_n)-\frac14I'(u_n)u_n}{\|u_n\|^2}=
\frac{a}{4}-\frac{\overline{\lambda}_n}{4}\|w_n\|_2^2
=\frac{a}{4}+o(1),
$$
we obtain $a=0$, this is a contradiction since $a>0$, therefore $\{u_n\}$ is
bounded in $H^1_0(\Omega)$.
\end{prooff}

\begin{lemma}\label{lem09}
Suppose that $a>0$, $b>0$ and $\lambda>a\lambda_1$.  Let $0<\delta<
b\mathcal{S}^2$. Then Eq.~\eqref{eq02} has a positive ground state
solution $u\in H^1_0(\Omega)$ with $I(u)<0$.
\end{lemma}
\begin{prooff}
Let $\varphi_{11}$ be the positive eigenfunction corresponding to
the first eigenvalue $\lambda_1$, then we see that
$$
a\|\varphi_{11}\|^2-\lambda\|\varphi_{11}\|_2^2=
(a-\frac{\lambda}{\lambda_1})\|\varphi_{11}\|^2<0\quad\text{and}\quad
b\|\varphi_{11}\|^4-\delta\|\varphi_{11}\|_4^4>(b-\frac{\delta}{\mathcal{S}^2})\|\varphi_{11}\|^4\geq0
$$
since $\lambda>a\lambda_1$ and $0<\delta\leq b\mathcal{S}^2$.
Therefore,
$$\aligned
I(t\varphi_{11})=&\frac{a}{2}\|t\varphi_{11}\|^2+\frac{b}{4}\|t\varphi_{11}\|^4-\frac{\lambda}{2}\|t\varphi_{11}\|^2_2-\frac{\delta}{4}\|t\varphi_{11}\|^4_4\\
=&\frac{1}{2}(a-\frac{\lambda}{\lambda_1})\|\varphi_{11}\|^2t^2+\frac14(b\|\varphi_{11}\|^4-\delta\|\varphi_{11}\|_4^4)t^4<0
\endaligned
$$
for $t>0$ small, i.e., there exists $t_0>0$ small such that
$I(t_0\varphi_{11})<0$.  On the other hand,
$$\aligned
I(u)=&\frac{a}{2}\|u\|^2+\frac{b}{4}\|u\|^4-\frac{\lambda}{2}\|u\|^2_2-\frac{\delta}{4}\|u\|^4_4\\
>&(a-\frac{\lambda}{\lambda_1})\|u\|^2+(b-\frac{\delta}{\mathcal{S}^2})\|u\|^4\to +\infty
\endaligned
$$
as $\|u\|\to \infty$, which means that $I$ is coercive.  Let
$c=\inf\{I(u): u\in H^1_0(\Omega)\}$.  By the Ekeland variational
principle (cf. \cite[Theorem 2.4]{W96}), there exists $(PS)_c$
sequence $\{u_n\}$ of $I$. Clearly, $-\infty<c<0$ and $\{u_n\}$ is
bounded.  Then by Proposition~\ref{pro01} and \eqref{eq69}, there
exists $u\in H^1_0(\Omega)\setminus\{0\}$ such that $u_n\to u$ in
$H^1_0(\Omega)$ as $n\to\infty$, up to a subsequence.  Therefore,
$u$ is a ground state solution of Eq.~\eqref{eq02} satisfying
$I(u)=c<0$.  Noting that $I(|u|)=I(u)$, so that $|u|$ is also a
global minimizer of $I$ on $H^1_0(\Omega)$, in this way we have that
$|u|$ is a ground state solution, then $|u|>0$ because of the
maximum principle, therefore we obtain a positive ground state
solution.
\end{prooff}

\begin{lemma}\label{lem07}
Assume that $a>0$, $b>0$ and $0<\delta<b\mathcal{S}^2$ with
$\lambda>a\lambda_1$. Let $\{\delta_n\}$ be a decreasing sequence
satisfying $\delta_n\to0$ as $n\to \infty$ and let $u_n$ be the
positive ground state solution corresponding to $\delta_n$, then
$\{u_n\}$ is bounded.
\end{lemma}
\begin{prooff}
Arguing by contradiction, we assume  $\|u_n\|\to\infty$ as $n\to\infty$.  Set $w_n=u_n/\|u_n\|$,
then $\|w_n\|=1$ and there exists $w_0\in H^1_0(\Omega)$ such that $w_n\rightharpoonup w_0$ in $H^1_0(\Omega)$
and $w_n\to w_0$ in $L^p(\Omega)$ $(1<p<4)$ as $n\to\infty$. Clearly,
$\|w_0\|\leq 1$.  Denote $c_n:=I(u_n)$.  It follows from  $I'(u_n)u_n=0$ that
$$
\lim_{n\to\infty}\frac{c_n}{\|u_n\|^4}=\lim_{n\to\infty}\frac{I(u_n)-\frac12I'(u_n)u_n}{\|u_n\|^4}
=\lim_{n\to\infty}\frac{\delta_n}{4}\|w_n\|_4^4-\frac{b}{4}=-\frac{b}{4}<0.
$$
On the other hand,
$$
\lim_{n\to\infty}\frac{c_n}{\|u_n\|^2}=\lim_{n\to\infty}\frac{I(u_n)-\frac14I'(u_n)u_n}{\|u_n\|^2}
=\frac{a}{4}-\lim_{n\to\infty}\frac{\lambda}{4}\|w_n\|_2^2=\frac{a}{4}-\frac{\lambda}{4}\|w_0\|_2^2.
$$
Thus we get a contradiction since $\|u_n\|\to\infty$ as
$n\to\infty$.  Therefore, $\{u_n\}$ is bounded.
\end{prooff}

\begin{lemma}\label{lem08}
Under the assumptions of Lemma~\ref{lem07} and let $c_n:=I(u_n)$,
then $c_n\to d$ as $n\to\infty$, where $d:=\inf\{J(u): u\in
H^1_0(\Omega)\}<0$,
$J(u)=\dfrac{a}{2}\|u\|^2+\dfrac{b}{4}\|u\|^4-\dfrac{\lambda}{2}\|u\|^2_2$
is given by \eqref{eq66}.
\end{lemma}
\begin{prooff}
By Lemma~\ref{lem07}, we see that $\{u_n\}$ is bounded, then
$\{c_n\}$ is also bounded.  Clearly, $\{c_n\}$ is nondecreasing on
$n$, thus $c_n\to c_0$ for some $c_0\in \bbr$ as $n\to\infty$.
Obviously, the functional $J$ is coercive, weakly low
semi-continuous and
$$
J(t\varphi_{11})=\frac12(a-\frac{\lambda}{\lambda_{1}})\|\varphi_{11}\|^2t^2+\frac{b}{4}\|\varphi_{11}\|^4t^4<0
$$
for $t>0$ small since $\lambda>a\lambda_1$.  Therefore, there exists
$v\in H^1_0(\Omega)\setminus\{0\}$ such that
$$
d=J(v)=\inf\{J(u): u\in H^1_0(\Omega)\}<0.
$$
It follows from Lemma~\ref{lem09} that $c_n\leq d$ since $I(u)\leq J(u)$
for each $u\in H^1_0(\Omega)$ and $\delta_n>0$, hence we get
$c_0\leq d$.  If $c_0<d$, then there exists $n_0\in \bbn$ such that
for each $n>n_0$,
$$
d\leq
J(u_n)=I(u_n)+\frac{\delta_n}{4}\|u_n\|_4^4=c_n+\frac{\delta_n}{4}\|u_n\|_4^4\leq
c_0+\frac{d-c_0}{2}=\frac{c_0+d}{2}<d
$$
since $\delta_n\to 0$ and $\{u_n\}$ is bounded, this is a
contradiction.  Therefore, we must have $c_0=d$.
\end{prooff}

\noindent{\bf Proof of Theorem~\ref{thm01}:}

Suppose that $0<\delta\leq b\mathcal{S}^2$. If Eq.~\eqref{eq02} has a nontrivial solution $u$ for each
$\lambda\in (0, a\lambda_1]$ then
$$\aligned
0=I'(u)u=& a\|u\|^2+b\|u\|^4-\lambda\|u\|^2_2-\delta\|u\|_4^4 \\
>&
(a-\frac{\lambda}{\lambda_1})\|u\|^2+(b-\frac{\delta}{\mathcal{S}^2})\|u\|^4\geq 0,
\endaligned
$$
this is impossible, therefore, Eq.~\eqref{eq02} has no
nontrivial solution for $\lambda\in(0, a\lambda_1]$.

Now we assume that $a>0$, $b>0$, $\lambda>a\lambda_1$ and $0<\delta<
b\mathcal{S}^2$.  Then it follows from Lemma~\ref{lem09} that
Eq.~\eqref{eq02} has a positive ground state solution.

Assume that $\lambda\in (a\lambda_k, a\lambda_{k+1}]$ with $k\geq
1$.  Let $\{\delta_n\}$ be a decreasing sequence satisfying
$\delta_n\to 0$ as $n\to \infty$. Let $u_n$ be the positive ground
state solution corresponding to $\delta_n$. Lemmas~\ref{lem07} and
\ref{lem08} show that
$$
J(u_n)=I(u_n)+\frac{\delta_n}{4}\|u_n\|^4_4=c_n+\frac{\delta_n}{4}\|u_n\|^4_4\to
c_0=d
$$
as $n\to\infty$, i.e., $\{u_n\}$ is a minimizing sequence of $J$.
Then there exists $0\leq\widetilde{u}_0\in H^1_0(\Omega)$ such that
$u_n\rightharpoonup \widetilde{u}_0$ in $H^1_0(\Omega)$ and $u_n\to
\widetilde{u}_0$ in $L^p(\Omega)$ $(1<p<4)$ as $n\to\infty$.
Therefore, we have
$$
d\leq J(\widetilde{u}_0)\leq \lim_{n\to\infty}J(u_n)=d,
$$
i.e., $\widetilde{u}_0$ is a global minimum of $J$ and then
$\widetilde{u}_0$ is a positive ground state solution of
\eqref{eq49} because of the maximum principle. It follows from Step
1 in the proof of Theorem~\ref{thm05} that
$\widetilde{u}_0=\overline{u}_1$, where $\overline{u}_1$ is given by
\eqref{eq65} with $j=1$. \qquad
\raisebox{-0.5mm}{\rule{1.5mm}{4mm}}\vspace{6pt}

\noindent{\bf Proof of Theorem~\ref{thm03}:}

\noindent(i)\ \ We set $V=H_1$ and $W=\overline{W}_{\varepsilon}$
(see \eqref{eq39}) with $\varepsilon$ so small that \eqref{eq63}
holds. Clearly, the assumptions $(a_1)$ and $(a_{33})$ of
Theorem~\ref{thm02} are satisfied.  It follows from
Lemmas~\ref{lem01} and \ref{lem06} that $(a_2)$ and $(a_{31})$ hold
with $\beta=(a\mathcal{S})^2/(4(\delta-b\mathcal{S}^2))$. If $u\in
V$ such that $\|u\|=\rho$ with $\rho>0$ small enough, then
$$
\aligned
I(u)=&\frac{a}{2}\|u\|^2+\frac{b}{4}\|u\|^4-\frac{\lambda}{2}\|u\|_2^2-\frac{\delta}{4}\|u\|_4^4\\
\geq &
\frac{a}{2}(1-\frac{\lambda}{a\lambda^+})\|u\|^2-\frac{\delta}{4\mathcal{S}^2}\|u\|^4\\
\geq& \kappa_0>0
\endaligned
$$
for some positive constant $\kappa_0$, that is, the
assumption $(a_{32})$ of Theorem~\ref{thm02} is also satisfied.
Since dim$(V\cap W)=1$ and $V\oplus W=H^1_0(\Omega)$,  by
Theorem~\ref{thm02}, we get that Eq.~\eqref{eq02} has one pair of
nontrivial solutions with the corresponding functional energy in
$[\kappa_0, (a\mathcal{S})^2/(4(\delta-b\mathcal{S}^2)))$.

We claim that the solutions  of Eq.~\eqref{eq02} obtained by the above are all sign-changing.

Indeed, otherwise, for the case of (1),  we may assume that there
exist $\{b_n\}$, $\{\overline{\lambda}_n\}$ with $b_n>0,
\overline{\lambda}_n\geq a\lambda_1$, $b_n\to 0,
\overline{\lambda}_n\to a\lambda_1$ as $n\to\infty$ and the
corresponding solutions $\{u_n\}$ with  $u_n\geq 0$ for all $n\in
\bbn$. By using Lemma~\ref{lem10}, we see that $\{u_n\}$ is bounded.
Recalling that the energy functional of \eqref{eq48} is defined as
$$
I_0(u)=\frac{a}{2}\|u\|^2-\frac{a\lambda_1}{2}\|u\|^2_2-\frac{\delta}{4}\|u\|^4_4,
\quad u\in H^1_0(\Omega),
$$
we deduce from Theorem~\ref{thm02} and Remark~\ref{rmk03}  that for $n$
large enough,
$$
0<\frac{\kappa_0}{2}\leq
I_0(u_n)=I(u_n)-\frac{b_n}{4}\|u_n\|^4+\frac{\overline{\lambda}_n-a\lambda_1}{2}\|u_n\|^2_2\leq\frac{(a\mathcal{S})^2}{4\delta}-\sigma
$$
which, together with Lemma~\ref{lem10}, implies that $I_0(u_n)\to
c_0\in[\kappa_0/2,(a\mathcal{S})^2/(4\delta)$, up to a subsequence.
On the other hand, for each $v\in H^1_0(\Omega)$,
$$
I'_0(u_n)v=I'(u_n)v-b_n\|u_n\|^2\int_{\Omega}\nabla u_n\nabla vdx -
(\overline{\lambda}_n-a\lambda_1)\int_{\Omega}u_nvdx=o(1),
$$
$$
I'_0(u_n)u_n=I'(u_n)u_n-b_n\|u_n\|^4 -
(\overline{\lambda}_n-a\lambda_1)\int_{\Omega}u^2_ndx=o(1),
$$
therefore, by using the standard arguments in \cite{BN83,S84,W96},
we obtain that there exists $u_0\in H^1_0(\Omega)\setminus\{0\}$
such that $u_n\to u_0$ in $H^1_0(\Omega)$ as $n\to\infty$, up to a
subsequence.

Let $\varphi_{11}$ be the positive eigenfunction of the first
eigenvalue $\lambda_1$, then we have
$$\aligned
0=&(a+b_n\|u_n\|^2)\int_{\Omega}\nabla u_n\nabla
\varphi_{11}dx-\overline{\lambda}_n\int_{\Omega}u_n\varphi_{11}dx-\delta\int_{\Omega}|u_n|^2u_n\varphi_{11}dx
\\
=&b_n\lambda_1\|u_n\|^2\int_{\Omega}u_n
\varphi_{11}dx-(\overline{\lambda}_n-a\lambda_{1})\int_{\Omega}u_n\varphi_{11}dx-\delta\int_{\Omega}|u_n|^2u_n\varphi_{11}dx\\
=&(b_n\lambda_1\|u_n\|^2-(\overline{\lambda}_n-a\lambda_1))\int_{\Omega}|u_n|\varphi_{11}dx-\delta\int_{\Omega}|u_n|^3\varphi_{11}dx\\
\leq& o(1)-\delta\int_{\Omega}|u_0|^3\varphi_{11}dx<0,
\endaligned
$$
this is a contradiction, so that our claim is true for the case of (1).  For the case (2), we observe that if $u$ is a
solution of \eqref{eq02}, then $v=\delta^{1/2}u$ is a solution of
$$
\left\{\aligned-(a+\frac{b}{\delta}\int_{\Omega}|\nabla v|^2dx)\Delta v &= \lambda v + |v|^{2}v, &\quad \text{in}\ \Omega, \\
v&=0,& \text{on}\ \partial\Omega,
\endaligned
\right.
$$
then by using similar arguments as above, we can show that $v$ is
sign-changing if $\delta$ large enough, so that $u$ is also
sign-changing.\vspace{6pt}

\noindent(ii)\ \ We now prove the asymptotical behavior of the
solutions as $b\to 0$.  Fix $a, \delta>0$.  Let $\{b_n\}$ and
$\{\overline{\lambda}_n\}$ be two sequences satisfying $b_n\to 0$
and $\overline{\lambda}_n\to a\lambda_1$ as $n\to \infty$ with
$\overline{\lambda}_n\geq a\lambda_1$.   Let $u_n$ be the solution
corresponding to $b_n$ and $\overline{\lambda}_n$ obtained above, it
follows from Lemma~\ref{lem10} that $\{u_n\}$ is bounded.  Similarly as in the proof of the above claim, we see that $u_n\to u_0\neq 0$ in
$H^1_0(\Omega)$ as $n\to\infty$ and $I_0'(u_0)v=0$ for all $v\in
H^1_0(\Omega)$, that is, $u_0$ is a nontrivial solution of
Eq.~\eqref{eq48}. A standard argument (cf. \cite{BN83}) shows that
$u_0$ is a sign-changing solution of Eq.~\eqref{eq48}. $\qquad
\raisebox{-0.5mm}{\rule{1.5mm}{4mm}}$

\section{Proof of Theorems \ref{thm04} and \ref{thm06}}\label{Sec04}
In this section, we are investigating the existence of solutions of
Eq.~\eqref{eq01} for $\lambda>a\lambda_1$ and $2<q<4$.  We regard
Eq.~\eqref{eq01} as a perturbed equation of the following equation
\begin{equation}\label{eq52}
\left\{\aligned-a\Delta u &= \lambda|u|^{q-2}u + \delta |u|^{2}u, &\quad \text{in}\ \Omega, \\
u&=0,& \text{on}\ \partial\Omega,
\endaligned
\right.
\end{equation}
then the corresponding functionals for Eq~\eqref{eq01} and Eq.~\eqref{eq52} can be respectively written as
$$
\overline{I}_b(u)=\frac{a}{2}\|u\|^2+\frac{b}{4}\|u\|^4-\frac{\lambda}{q}\int_{\Omega}|u|^qdx-\frac{\delta}{4}\int_{\Omega}|u|^4dx,
$$
and
$$
\overline{I}_0(u)=\frac{a}{2}\|u\|^2-\frac{\lambda}{q}\int_{\Omega}|u|^qdx-\frac{\delta}{4}\int_{\Omega}|u|^4dx.
$$
It follows from \cite{BN83} (see also \cite{S84}) that
Eq.~\eqref{eq52} has a positive ground state solution $u_0\in
H^1_0(\Omega)$, thus we obtain
$$
a\|u_0\|^2=\lambda\int_{\Omega}|u_0|^qdx+\delta\int_{\Omega}|u_0|^4dx
$$
and
\begin{equation}\label{eq53}
\overline{I}_0(u_0)=\frac{q-2}{2q}a\|u_0\|^2+\frac{4-q}{4q}\delta\int_{\Omega}|u_0|^4dx.
\end{equation}
Moreover, $u_0$ is a mountain pass solution and $\overline{I}_0$ possesses the
following properties:
\begin{itemize}
\item [(M1)] there exist $c$, $r>0$ such that if $\|u\|=r$  then $\overline{I}_0(u)\geq
c$ and there exists $v_0\in H^1_0(\Omega)$ such that $\|v_0\|>r$ and
$\overline{I}_0(v_0)<0$;
\item [(M2)] there exists a critical point $u_0\in H^1_0(\Omega)$ of
$\overline{I}_0$ such that
$$
\overline{I}_0(u_0)=c_0:=\min_{\gamma\in \Gamma}\max_{t\in
[0,1]}\overline{I}_0(\gamma(t)),
$$
where $\Gamma:=\{\gamma\in C([0,1], H^1_0(\Omega)): \gamma(0)=0,
\gamma(1)=v_0\}$;
\item [(M3)] $0<c_0:=\inf\{\overline{I}_0(u): \overline{I}_0'(u)=0, u\in
H^1_0(\Omega)\setminus\{0\}\}<(a\mathcal{S})^2/(4\delta)$;
\item [(M4)] the set $\overline{S}:=\{u\in H^1_0(\Omega): \overline{I}_0'(u)=0,
\overline{I}_0(u)=c_0\}$ is compact in $H^1_0(\Omega)$;
\item [(M5)] there exists a path $\gamma_0\in \Gamma$ passing
through $u_0$ at $t=t_0$ and satisfying
$$
\overline{I}_0(u_0)>\overline{I}_0(\gamma_0(t))\ \ \text{for all}\
t\neq t_0.
$$
\end{itemize}
In fact, we can take $v_0=Tu_0$ with $T>\sqrt{q/2}$ in $(M1)$, and
set $\gamma_0(t)=tv_0$ with $t_0=1/T$ in $(M5)$.  It is easy to
verify that $(M2)$ and $(M3)$ hold.  To check $(M4)$, let $\{u_n\}\subset \overline{S}$ be a sequence, we deduce  from \eqref{eq53} and the definition of $\overline{S}$ that
$\{u_n\}$ is a bounded (PS) sequence of
$\overline{I}_0$ at the level $c_0$.  Noting that $c_0<(a\mathcal{S})^2/(4\delta)$, by
using the arguments of \cite{BN83,S84}, we see that
$\overline{I}_0$ satisfies the Palais-Smile condition, and then
there exists $\overline{u}_0\in H^1_0(\Omega)\setminus\{0\}$ such
that $u_n\to \overline{u}_0$ in $H^1_0(\Omega)$ as $n\to\infty$, so
that $\overline{I}_0(\overline{u}_0)=c_0$ and $\overline{u}_0\in
\overline{S}$.\vspace{6pt}

A nature way to construct a (PS)$_c$ sequence of $\overline{I}_b$ with $c<(a\mathcal{S})^2/(4(\delta-b\mathcal{S}^2))$ can be processed as follows: let $v_0\in H_0^1(\Omega)$ be as in (M1), since $\overline{I}_b(v_0)<0$ for $b>0$ small enough, it is easy to
show that
$\overline{I}_b$ has the mountain pass geometry and hence there exists a
sequence $\{u_n\}\subset H_0^1(\Omega)$ satisfying $\overline{I}_b(u_n)\to
\overline{c}_{b}$ and $\overline{I}_b'(u_n)\to 0$ in
$H^{-1}(\Omega)$, where
$$
\overline{c}_b=\min_{\gamma\in \Gamma}\max_{t\in
[0,1]}\overline{I}_b(\gamma(t))
$$
and $\Gamma$ is defined in (M2). Let $b>0$ be so small that we have the following estimate of $\overline{c}_b$:
$$\aligned
\overline{c}_b\leq \max_{t\in [0,T]}\overline{I}_b(tu_0)\leq&
\max_{t\in[0, T]}\overline{I}_0(tu_0)+\frac{bT^4}{4}\|u_0\|^4 \\
=& c_0+\frac{bT^4}{4}\|u_0\|^4 \\
<&\frac{(a\mathcal{S})^2}{4\delta}
<\frac{(a\mathcal{S})^2}{4(\delta-b\mathcal{S}^2)}.
\endaligned
$$
However, we are not sure whether the (PS)$_{\overline{c}_b}$ sequence of the functional $\overline{I}_b$ has a convergent subsequence in $H_0^1(\Omega)$ or not because the exponent of the nonlocal term $\|u\|^4$ is equal
to the critical Sobolev exponent when $N=4$, which leads to the
difficulty of proving the boundedness of $\{u_n\}$.

We borrow ideas of \cite{BJ07,CZ14,CZZ15,JS14} and define a modified
mountain pass level of $\overline{I}_b$ by
$$
c_b:=\min_{\gamma\in \Gamma_M}\max_{t\in
[0,1]}\overline{I}_b(\gamma(t))
$$
where
\begin{equation}\label{eq4.4}
\Gamma_M:=\{\gamma\in \Gamma: \sup_{t\in
[0,1]}\|\gamma(t)\|\leq M\}\ \text{with}\ M=2T\|u_0\|^2.
\end{equation}
It follows from
the choice of $M$ that $\gamma_0\in \Gamma_M$ and
$c_0=\min_{\gamma\in\Gamma_M}\max_{t\in
[0,1]}\overline{I}_0(\gamma(t))$, where $c_0$ and $\gamma_0$ are respectively given in (M3) and (M5).  Moreover, we have the following result:
\begin{lemma}\label{lem12}
The mountain pass level $c_b$ is continuous at $0$, i.e.,
$\lim_{b\to0}c_b=c_0$.
\end{lemma}
\begin{prooff}
Clearly, $c_b\geq c_0$ holds.  It follows from $(M5)$ that
$$
\aligned c_b\leq\max_{1\leq
t\leq1}\overline{I}_b(tv_0)\leq&\max_{1\leq t\leq
T}\overline{I}_b(tu_0)+\frac{bT^4}{4}\|u_0\|^4 \\
=& c_0+\frac{bT^4}{4}\|u_0\|^4 \\
=&c_0+o(1)
\endaligned
$$
as $b\to 0$.
\end{prooff}

Given $d>0$, a set $A\subset H^1_0(\Omega)$ and a point $u\in H^1_0(\Omega)$, we denote
$$
B_d(u):=\{v\in H^1_0(\Omega): \|v-u\|\leq d\},\quad
A^d:=\bigcup_{u\in A}B_d(u).
$$
\begin{lemma}\label{lem13}Let $\overline{S}$ be given in $(M4)$.
For fixed $d>0$, if a sequence $\{u_j\}\subset \overline{S}^d$ then
$\{u_j\}$ converges weakly to some $u\in \overline{S}^{2d}$ as $j\to\infty$, up to a
subsequence.
\end{lemma}
\begin{prooff}
By the definition of $\overline{S}^{d}$, there
exists $\{v_j\}\subset \overline{S}$ such that $\|u_j-v_j\|\leq d$ for all $j\in\mathbb{N}$.  It follows from the compactness that
there exists $v\in \overline{S}$ such that $\{v_j\}$ converges strongly to $v$ in $H^1_0(\Omega)$ as $j\to\infty$, up to a
subsequence, so that $\|u_j-v\|\leq 2d$ for large $j$, that is, $\{u_j\}\subset B_{2d}(v)$ for $j$ large enough.  Since $B_{2d}(v)$ is weakly closed in
$H^1_0(\Omega)$ and $\{u_j\}$ is clearly bounded in $H_0^1(\Omega)$, we get that there exists $u\in \overline{S}^{2d}$ such that $u_j\rightharpoonup u$  as $j\to\infty$, up to a
subsequence.
\end{prooff}
\begin{remark}\label{rmk06}{\em By the definition of $\overline{S}$ in (M4), we see that there exists a positive constant $C$ such that $\|u\|\geq C$ for all $u\in \overline{S}$ since, otherwise, we must have a sequence $\{u_n\}\subset \overline{S}$ such that
$u_n\to 0$ in $H_0^1(\Omega)$ as $n\to\infty$, which yields $c_0=0$, a contradiction with the definition of $c_0$. Therefore, we may choose a $d>0$ so small that $u\neq 0$ for all $u\in \overline{S}^{2d}$.
}
\end{remark}
\begin{lemma}\label{lem14}
Assume that $\{b_j\}$ is a sequence of positive numbers with $b_j\to
0$ as $j\to\infty$.  Let $d>0$ be small enough and let $\{u_j\}\subset \overline{S}^{d}$ be such that
\begin{equation}\label{eq55}
\lim_{j\to\infty}\overline{I}_{b_j}(u_j)\leq c_0, \quad
\lim_{j\to\infty}\overline{I}'_{b_j}(u_j)=0.
\end{equation}
Then $\{u_j\}$
strongly converges to some $u\in \overline{S}$, up to a
subsequence.
\end{lemma}
\begin{prooff}
By Lemma~\ref{lem13} and Remak~\ref{rmk06}, up to a subsequence, we get that there exists $u\in \overline{S}^{2d}$ such that
$u_j\rightharpoonup u$ in $H_0^1(\Omega)$ as $j\to\infty$ and $u\neq 0$.  It can be deduce from \eqref{eq55} and the boundedness of
$\{u_j\}$  that for all $v\in H^1_0(\Omega)$,
$$\aligned
\overline{I}'_0(u)v=&a\int_{\Omega}\nabla u\nabla vdx -
\lambda\int_{\Omega}|u|^{q-2}uvdx-\delta\int_{\Omega}|u|^2uvdx\\
=&\lim_{j\to\infty}a\int_{\Omega}\nabla u_j\nabla vdx -
\lambda\int_{\Omega}|u_j|^{q-2}u_jvdx-\delta\int_{\Omega}|u_j|^2u_jvdx\\
=&\lim_{j\to\infty}\overline{I}'_{b_j}(u_j)v-b_j\|u_j\|^2\int_{\Omega}\nabla
u_j\nabla vdx=0,
\endaligned
$$
hence $\overline{I}'_0(u)=0$.  Furthermore, since $\{u_j\}\subset
\overline{S}^{d}$, we obtain that
$$
\overline{I}'_0(u_j)=\overline{I}'_{b_j}(u_j)-b_j\|u_j\|^2 u_j=o(1).
$$
On the other hand,
\begin{equation}\label{eq56}
c_0\geq
\lim_{j\to\infty}\overline{I}_{b_j}(u_j)=\lim_{j\to\infty}\overline{I}_{0}(u_j)+\lim_{j\to\infty}b_j\|u_j\|^4
=\lim_{j\to\infty}\overline{I}_{0}(u_j):= m,
\end{equation}
then $\{u_j\}$ is a $(PS)_m$ sequence of $\overline{I}_0$.
Therefore, up to a subsequence,
$$\aligned
\overline{I}_{0}(u)=&\frac{a}{2}\int_{\Omega}|\nabla
u|^2dx-\frac{\lambda}{q}\int_{\Omega}|u|^qdx-\frac{\delta}{4}\int_{\Omega}|u|^4dx\\
=&\frac{q-2}{2q}a\int_{\Omega}|\nabla
u|^2dx+\frac{4-q}{4q}\delta\int_{\Omega}|u|^4dx\\
\leq&\lim_{j\to\infty}\frac{q-2}{2q}a\int_{\Omega}|\nabla
u_j|^2dx+\frac{4-q}{4q}\delta\int_{\Omega}|u_j|^4dx\\
=&\lim_{j\to\infty}(\overline{I}_{0}(u_j)-\frac1q\overline{I}'_0(u_j)u_j)\\
=&m.
\endaligned
$$
It follows from $(M3)$ and \eqref{eq56} that
$m=\overline{I}_0(u)=c_0$, which implies that $u\in
\overline{S}$.  A standard argument (cf. \cite{BN83,S84})
shows that $u_j\to u$ in $H^1_0(\Omega)$ as $j\to\infty$.
\end{prooff}

Now we set $D_b:= \max\{\overline{I}_b(\gamma_0(t)): 0\leq t\leq
1\}$.  Clearly, $c_b\leq D_b$ and  moreover, it holds $\lim_{b\to 0}D_b\leq
c_0$, so that, by Lemma~\ref{lem12}, we have
\begin{equation}\label{eq58}
\lim_{b\to 0}c_b=\lim_{b\to 0}D_b=c_0.
\end{equation}
The following Lemmas~\ref{lem15}-\ref{lem17}
are quite analogous to Propositions 3-5 in \cite{JS14}. However,
some arguments in \cite{JS14} are based on the assumptions that the
perturbation operator is compact, which is not our case here,
therefore we prefer to give the modified proofs for sake of
completeness.

\begin{lemma}\label{lem15}
Let $d_1>d_2>0$ be two small constants, then there exist positive numbers $\alpha$
and $b_0$ depending on $d_1$ and $d_2$ such that for each $b\in(0,
b_0)$, the following is true:
$$
\|\overline{I}'_b(u)\|\geq \alpha \ \quad \text{for all}\  u\in
\overline{I}^{D_b}_b\cap(\overline{S}^{d_1}\setminus\overline{S}^{d_2}),
$$
where
$\overline{I}_b^{D_b}=\{u\in H^1_0(\Omega): \overline{I}_b(u)\leq
D_b\}$ as usual.
\end{lemma}
\begin{prooff}
Argue by contradiction, we may suppose that there exist  a sequence
$\{b_j\}$ of positive numbers with $\lim_{j\to\infty}b_j=0$ and a sequence of functions
$\{u_j\}\subset
\overline{I}^{D_{b_j}}_{b_j}\cap(\overline{S}^{d_1}\setminus\overline{S}^{d_2})$
such that $\lim_{j\to\infty}\overline{I}'_{b_j}(u_j)=0$. Then by
\eqref{eq58}, we obtain that $\lim_{j\to\infty}\overline{I}_{b_j}(u_j)\leq c_0$.
It follow from Lemma~\ref{lem14}  that there exists $u\in
\overline{S}$ such that $u_j\to u$ in $H^1_0(\Omega)$. As
a consequence, dist$(u_j, \overline{S})\to 0$ as $j\to
\infty$, this contradict $u_j\notin \overline{S}^{d_2}$.
\end{prooff}

\begin{lemma}\label{lem16}
Let $d>0$ be fixed. Then there exists $\delta>0$ such
that if $b>0$ small enough,
$$
t\in[0, 1]\ \text{and}\ \overline{I}_b(\gamma_0(t))\geq c_b-\delta
\Rightarrow \gamma_0(t)\in \overline{S}^d.
$$
\end{lemma}
\begin{prooff}
Since the proof is quite similar to that of \cite{JS14}.  We omit
the details.
\end{prooff}

\begin{lemma}\label{lem17}
Let $d>0$ be a small constant.  Then there exist $b>0$ sufficiently
small, depending on $d$, and a sequence $\{u_j\}\subset
\overline{S}^d\cap\overline{I}^{D_b}_b$ such that
$\overline{I}'_b(u_j)\to 0$ as $j\to \infty$.
\end{lemma}
\begin{prooff}
Arguing by contradiction, we may assume that there exist two
sequences $\{b_j\}$ and $\{c_j\}$ of positive numbers with $b_j\to
0$ such that $\|\overline{I}'_{b_j}(u)\|\geq c_j>0$ for all $u\in
\overline{S}^d\cap\overline{I}^{D_{b_j}}_{b_j}$.  By
Lemmas~\ref{lem14} and \ref{lem15}, there exists some $\alpha>0$,
independent of $j$, such that
$$
\|\overline{I}'_{b_j}(u)\|\geq \alpha\  \text{for all}\ u\in
\overline{I}_{b_j}^{D_{b_j}}\cap(\overline{S}^d\setminus\overline{S}^{\frac{d}{2}}),
$$
which gives that $\{c_j\}$ has a positive lower bound.  Moreover,
there exists $k>0$, independent of $j$, such that
$\|\overline{I}'_{b_j}(u)\|\leq k$ for all $u\in
\overline{S}^d$ since $\overline{S}^d$ is bounded in $H_0^1(\Omega)$.  By using
Lemma~\ref{lem16} and \eqref{eq58}, we can choose $\delta>0$ so small and $j$ so large that
\begin{equation}\label{eq57}
t\in[0,1], \overline{I}_{b_j}(\gamma_0(t))\geq c_b-\frac{\delta}{4}\
\text{implies}\ \gamma_0(t)\in \overline{S}^{\frac{d}{2}}
\end{equation}
and
\begin{equation}\label{eq59}
D_{b_j}-c_{b_j}<\frac{\delta}{4}\ \text{and}\
D_{b_j}-c_{b_j}<\frac{\alpha^2d}{4k}-\frac{\delta}{4}.
\end{equation}

From now on, we fix a $j$ so large that \eqref{eq57} and \eqref{eq59}
are true, and we denote $b_j$ just by $b$.
Now, consider a pseudo-gradient vector field $V_b$ of
$\overline{I}_b$ and a neighborhood $\mathcal{N}_b$ of
$\overline{S}^d\cap\overline{I}^{D_b}_b$ satisfying
$\mathcal{N}_b\subset B_{M}(0)$, where $M$ is given by \eqref{eq4.4}.  We observe that $D_b<M$ for $b$
small enough.

Let us respectively define two functions $\eta_b\in C^{1,1}(H^1_0(\Omega), [0, 1])$  and $\xi_b\in C^{1,1}(\bbr, [0, 1])$ by
$$
\eta_b=\left\{\aligned &1, &\quad \text{on}\ \overline{S}^d\cap\overline{I}^{D_b}_b, \\
&0,& \text{on}\ H^1_0(\Omega)\setminus\mathcal{N}_b
\endaligned
\right.
$$
and
$$
\xi_b(t)=\left\{\aligned &1, &\quad \text{if}\ |t-c_b|\leq\frac{\delta}{2}, \\
&0,& \text{if}\ |t-c_b|\geq\delta.
\endaligned
\right.
$$
Then it is clear that the following  Cauchy initial value problem
$$
\left\{\aligned\frac{d}{d t}\psi_b(u,t)&=-\eta_b(\psi_b(u, t))\xi_b(\overline{I}_b(\psi_b(u,t)))V_b(\psi_b(u,t)), \\
\psi_b(u,0)&=u
\endaligned
\right.
$$
has a global solution $\psi_b:
H^1_0(\Omega)\times \bbr\to H^1_0(\Omega)$ since  $\eta_b$ and $\xi_b\circ \overline{I}_b\in
C^{1,1}(H^1_0(\Omega), [0,1])$ and $V_b$ is locally Lipschitz
continuous.

Now we claim that for all $t\in[0, 1]$,
$\overline{I}_b(\psi_b(\gamma_0(t),t_b))\leq c_b-\delta/4$, where
$t_b:= \delta/(2c^2)$.  Indeed, if for each  $t\in[0, 1]$, there
exists $t_0\leq t_b$ such that $\overline{I}_b(\psi_b(\gamma_0(t),t_0))\leq c_b-\delta/4$,
then we get the claim since it follows from
\begin{equation}\label{eq60}
\frac{d}{d\tau}\overline{I}_b(\psi_b(\gamma_0(t),\tau))\leq
-\eta_b(\psi_b(\gamma_0(t),\tau))\xi_b(\overline{I}_b(\psi_b(\gamma_0(t),\tau)))\|\overline{I}'_b(\psi_b(\gamma_0(t),\tau))\|^2
\end{equation}
that
$$
\overline{I}_b(\psi_b(\gamma_0(t), t_b))\leq
\overline{I}_b(\psi_b(\gamma_0(t), t_0))\leq c_b-\frac{\delta}{4}.
$$
Otherwise, there exists some $t\in[0,1]$ such that
\begin{equation}\label{eq4.11}
\overline{I}_b(\psi_b(\gamma_0(t), \tau))>c_b-\frac{\delta}{4}\quad
\text{for all}\ \tau\in[0, t_b].
\end{equation}
Thus by \eqref{eq57}, we have
$$
\gamma_0(t)=\psi_b(\gamma_0(t), 0)\in
\overline{S}^{\frac{d}{2}}\ \text{and}\
\xi_b(\overline{I}_b(\psi_b(\gamma_0(t), \tau)))=1\ \text{for all}\
\tau\in[0, t_b].
$$
If $\psi_b(\gamma_0(t), \tau)\in \overline{S}^d$ ($\forall\,\tau
\in[0, t_b]$) then it follows from \eqref{eq59} and
$\overline{I}_b(\psi_b(\gamma_0(t), \tau))\leq
\overline{I}_b(\gamma_0(t))\leq D_b$ that $\eta_b(\psi_b(\gamma_0(t), \tau))=1$
($\forall\,\tau
\in[0, t_b]$), which, together with \eqref{eq60}, implies that
$$
\frac{d}{dt}\overline{I}_b(\psi_b(\gamma_0(t), \tau))\leq -c^2,
$$
where $c:=c_j$ for the above fixed $j$.  Thus we obtain
$$
\overline{I}_b(\psi_b(\gamma_0(t), t_b))=
\int_0^{t_b}\frac{d}{dt}\overline{I}_b(\psi_b(\gamma_0(t),
\tau))d\tau+ \overline{I}_b(\gamma_0(t)) \leq
D_b-\int_0^{t_b}c^2d\tau =D_b-t_bc^2 <c_b-\frac{\delta}{4},
$$
a contradiction with \eqref{eq4.11}.  Therefore we must have
$\psi_b(\gamma_0(t),
\tau_0)\notin \overline{S}^d$ for some $\tau_0\in [0, t_b]$, so that there exist $\tau_1,
\tau_2: 0\leq \tau_1<\tau_2\leq \tau_0$ satisfying
$\psi_b(\gamma_0(t), \tau_1)\in
\partial\overline{S}^{\frac{d}{2}}$, $\psi_b(\gamma_0(t), \tau_2)\in
\partial\overline{S}^{d}$ and $\psi_b(\gamma_0(t), \tau)\in
\overline{S}^d\setminus\overline{S}^{\frac{d}{2}}$
for all $\tau\in (\tau_1, \tau_2)$.  It is not hard to show that
$\tau_2-\tau_1\geq d/(4k)$ since, otherwise, we shall get, by
the definition of pseudo-gradient vector field, that
$$
\aligned \frac{d}{2}=&\|\psi_b(\gamma_0(t),
\tau_2)-\psi_b(\gamma_0(t), \tau_1)\|
=\|\int_{\tau_1}^{\tau_2}\frac{d}{d\tau}\psi_b(\gamma_0(t),
\tau)d\tau\|\\ \leq
&\int_{\tau_1}^{\tau_2}\|\frac{d}{d\tau}\psi_b(\gamma_0(t),
\tau)\|d\tau \leq \int_{\tau_1}^{\tau_2}\|I'_b(\psi_b(\gamma_0(t),
\tau))\|d\tau\leq \int_{\tau_1}^{\tau_2}2kd\tau < \frac{d}{2},
\endaligned
$$
this is impossible.  Thus, from \eqref{eq59} and
\eqref{eq60}, we obtain
$$
\aligned \overline{I}_b(\psi_b(\gamma_0(t), t_b))\leq &
\overline{I}_b(\psi_b(\gamma_0(t), \tau_2))
=\int_0^{\tau_2}\frac{d}{d\tau}\overline{I}_b(\psi_b(\gamma_0(t),
\tau))d\tau + \overline{I}_b(\gamma_0(t))\\ \leq &
\int_{\tau_1}^{\tau_2}\frac{d}{d\tau}\overline{I}_b(\psi_b(\gamma_0(t),
\tau))d\tau + \overline{I}_b(\gamma_0(t)) \leq
D_b-\alpha^2(\tau_2-\tau_1)\\ <&D_b-\frac{d\alpha^2}{4k}
<c_b-\frac{\delta}{4},
\endaligned
$$
a contradiction with \eqref{eq4.11} again. Hence we have prove the claim.

Finally, we set $\widetilde{\gamma}_0(t):=\psi_b(\gamma_0(t), t_b)$.
Then $\widetilde{\gamma}_0(t)\in \Gamma_M$ and
$\overline{I}_b(\widetilde{\gamma}_0(t))<c_b$ for all $t\in [0, 1]$,
this contradicts the definition of $c_b$. Thus we complete the proof of this lemma.
\end{prooff}

\noindent{\bf Proof of Theorem~\ref{thm04}:}

Let us fix $d>0$ small enough.  By Lemma~\ref{lem17}, \eqref{eq4.4}
and (M3), we deduce that there exist $b>0$ sufficiently small and a
(PS) sequence $\{u_n\}$ of $\overline{I}_b$ with $\{u_n\}\subset
\overline{S}^{\frac{d}{2}}\cap\overline{I}^{D_b}_b$ and
$D_b<(a\mathcal{S})^2/(4\delta)<(a\mathcal{S})^2/(4(\delta-b\mathcal{S}^2))$.
It follows from the \eqref{eq53}, $(M3)$ and the definition of
$\overline{S}^d$ that $\{u_n\}$ is bounded in $H^1_0(\Omega)$, in
particular, $\|u_n\|^2\leq 1+(aq\mathcal{S}^2)/(2\delta(q-2))$ for
all $n\in \bbn$.  Therefore, there exists $c>0$ with $
c<(a\mathcal{S})^2/(4\delta)<(a\mathcal{S})^2/(4(\delta-b\mathcal{S}^2))$
such that $\overline{I}_b(u_n)\to c$, up to a subsequence.  Next we will
show that there exists $\overline{u}_0\in \overline{S}^d$ such that
$u_n\to \overline{u}_0$ in $H^1_0(\Omega)$ as $n\to\infty$, up to a
subsequence, provided one of the three conditions $(C1)$-$(C3)$
holds.  Arguing by contradiction, if $\{u_n\}$ has no convergent
subsequence in $H^1_0(\Omega)$  then by Proposition~\ref{pro01},
there exist a weak limit $\overline{u}_0\in H^1_0(\Omega)$ of
$\{u_n\}$, a number $k\in \bbn$ and further, for every $i\in
\{1,2,\cdots,k\}$, a sequence of values $\{R^i_n\}\subset \bbr^+$,
points $\{x^i_n\}\subset\overline{\Omega}$ and a function $v_i\in
\mathcal{D}^{1,2}(\bbr^4)$ satisfying \eqref{eq07}, \eqref{eq08} and
$$
\overline{I}_b(u_n)=\widetilde{I}_b(\overline{u}_0)+\sum_{i=1}^k\widetilde{I}^{\infty}_b(v_i)+o(1)
$$
up to a subsequence, where $o(1)\to 0$ as $n\to \infty$,
$\widetilde{I}_b$ and $\widetilde{I}^{\infty}_b$ are respectively
given by \eqref{eq09} and \eqref{eq10}.  On one hand, by
\eqref{eq07} and \eqref{eq09}, we have
\begin{equation}\label{eq54}\aligned
\widetilde{I}_b(\overline{u}_0)=&\widetilde{I}_b(\overline{u}_0)-\frac14\Big((a+b(\|\overline{u}_0\|^2+
\sum^k_{i=1}\|v_i\|^2_{\mathcal{D}^{1,2}}))\|\overline{u}_0\|^2-\lambda\int_{\Omega}|\overline{u}_0|^qdx-\delta\int_{\Omega}|\overline{u}_0|^4dx\Big)\\
=&\frac{a}{4}\|\overline{u}_0\|^2-\frac{(4-q)\lambda}{4q}\int_{\Omega}|\overline{u}_0|^qdx\\
\geq&\frac{a}{4}\|\overline{u}_0\|^2-\frac{(4-q)\lambda|\Omega|^{\frac{4-q}{4}}}{4q\mathcal{S}^{q/2}}\|\overline{u}_0\|^q.
\endaligned
\end{equation}
On the other hand, it can be deduced from \eqref{eq10} that for
$j\in\{1,2,\cdots,k\}$,
$$\aligned
0=&(a+b(\|\overline{u}_0\|^2+\sum^k_{i=1}\|v_i\|^2_{\mathcal{D}}))\|v_j\|^2_{\mathcal{D}}-\delta\int_{\bbr^4}|v_j|^4dx\\
\geq&(a+b\|\overline{u}_0\|^2)\|v_j\|^2_{\mathcal{D}}-\frac{\delta-b\mathcal{S}^2}{\mathcal{S}^2}\|v_j\|^4_{\mathcal{D}},
\endaligned
$$
which implies that
$$
\|v_j\|^2_{\mathcal{D}^{1,2}}\geq\frac{\mathcal{S}^2(a+b\|\overline{u}_0\|^2)}{\delta-b\mathcal{S}^2}.
$$
Therefore, by \eqref{eq08}, we obtain that
$$\aligned
\widetilde{I}^{\infty}_b(v_j)=&\widetilde{I}^{\infty}_b(v_j)-\frac{1}{4}(a+b(\|\overline{u}_0\|^2
+\sum^k_{i=1}\|v_i\|^2_{\mathcal{D}})\|v_j\|^2_{\mathcal{D}}-\delta\int_{\bbr^4}|v_j|^4dx)\\
=&\frac{a}{4}\|v_j\|^2_{\mathcal{D}}
\geq\frac{a\mathcal{S}^2(a+b\|\overline{u}_0\|^2)}{4(\delta-b\mathcal{S}^2)}.
\endaligned
$$
this, together with \eqref{eq54}, yields that
$$\aligned
c=&\lim_{n\to\infty}\overline{I}_b(u_n)
\geq\widetilde{I}_b(\overline{u}_0)+\widetilde{I}^{\infty}_b(v_i)\\
=&\frac{a}{4}\|\overline{u}_0\|^2-\frac{(4-q)\lambda}{4q\mathcal{S}_q^{q/2}}\|\overline{u}_0\|^q+\frac{(a\mathcal{S})^2}{4(\delta-b\mathcal{S})^2}
+\frac{ab\mathcal{S}^2}{4(\delta-b\mathcal{S}^2)}\|\overline{u}_0\|^2\\
=&\frac{(a\mathcal{S})^2}{4(\delta-b\mathcal{S})^2}+\Big(\frac{a\delta}{4(\delta-b\mathcal{S}^2)}
-\frac{(4-q)\lambda}{4q\mathcal{S}_q^{q/2}}\|\overline{u}_0\|^{q-2}\Big)\|\overline{u}_0\|^2\\
\geq&
\frac{(a\mathcal{S})^2}{4(\delta-b\mathcal{S})^2}+\Big(\frac{a\delta}{4(\delta-b\mathcal{S}^2)}
-\frac{(4-q)\lambda|\Omega|^{\frac{4-q}{4}}}{4q\mathcal{S}^{q/2}}\Big(1+\frac{aq\mathcal{S}^2}{2\delta(q-2)}\Big)^{\frac{q-2}{2}}\Big)\|\overline{u}_0\|^2,
\endaligned
$$
here we have used the inequality $\|\overline{u}_0\|^2\leq
1+(aq\mathcal{S}^2)/(2\delta(q-2))$ since $\|u_n\|^2\leq
1+(aq\mathcal{S}^2)/(2\delta(q-2))$ and $\overline{u}_0$ is the weak
limit of $u_n$ in $H_0^1(\Omega)$.

Now, recalling that $2<q<4$, it follows from one of the assumptions
$(C1)$-$(C3)$ that
$$
c\geq \frac{(a\mathcal{S})^2}{4(\delta-b\mathcal{S})^2},
$$
this is a contradiction since $c<(a\mathcal{S})^2/(4(\delta-b\mathcal{S}^2))$.
Therefore, up to a subsequence, $\{u_n\}$ converges to some
$\overline{u}_0\in H^1_0(\Omega)$.  Noting that
$\overline{u}_0\neq 0$ (cf. Remak~\ref{rmk06}), we know that
$\overline{u}_0$ is a desired solution of Eq.~\eqref{eq01}.
$\qquad \raisebox{-0.5mm}{\rule{1.5mm}{4mm}}$\vspace{6pt}

\noindent{\bf Proof of Theorem~\ref{thm06}:}

\noindent $(i)$\quad We suppose on the contrary that Eq.~\eqref{eq01} has a solution $u\in
H_0^1(\Omega)\setminus\{0\}$ under assumptions $\delta<b\mathcal{S}^2$ and \eqref{eq1.5}, that is,
$\overline{I}'_b(u)v=0$ for
all $v\in H^1_0(\Omega)$, so that
\begin{equation}\label{eq68}
\overline{I}'_b(u)u=a\|u\|^2+b\|u\|^4-\lambda\|u\|_q^q-\delta\|u\|_4^4=0.
\end{equation}
For $t>0$, we set
$f(t)=\overline{I}_b(tu)$, then
$$
f'(t)=at\|u\|^2+t^3(b\|u\|^4-\delta\|u\|_4^4)-\lambda
t^{q-1}\|u\|_q^q:=th(t),
$$
where
$h(t)=a\|u\|^2+t^2(b\|u\|^4-\delta\|u\|_4^4)-\lambda
t^{q-2}\|u\|_q^q$.  It follows from \eqref{eq68} that $f'(1)=0$.  On the other hand, it is easy to see that
$$
\overline{t}=\bigg(\frac{\lambda(q-2)\|u\|_q^q}{2(b\|u\|^4-\delta\|u\|_4^4)}\bigg)^{\frac1{4-q}}
$$
is the unique root of equation $h'(t)=0$ on $(0,+\infty)$, which, together with the H\"{o}lder inequality, the inequality \eqref{eq1.5}  and the fact of $q\in (2,4)$, implies
$$
\aligned
\min_{t>0}h(t)=h(\overline{t})=&a\|u\|^2+\overline{t}^2(b\|u\|^4-\delta\|u\|_4^4)-\lambda
\overline{t}^{q-2}\|u\|_q^q \\
=&a\|u\|^2-\frac{(4-q)}{2}\lambda\|u\|_q^q\Big(\frac{\lambda(q-2)\|u\|_q^q}{2(b\|u\|^4-\delta\|u\|_4^4)}\Big)^{\frac{q-2}{4-q}}\\
\geq&a\|u\|^2-\frac{(4-q)}{2}\Big(\frac{(q-2)\mathcal{S}^2}{2(b\mathcal{S}^2-\delta)}\Big)^{\frac{q-2}{4-q}}
\frac{\lambda^{\frac{2}{4-q}}\|u\|_q^{\frac{2q}{4-q}}}{\|u\|^{\frac{4(q-2)}{4-q}}}\\
\geq&\|u\|^2\Big(a-(4-q)\Big(\frac{\lambda}{2}\Big)^{\frac{2}{4-q}}\Big(\frac{|\Omega|^{\frac{4-q}{4}}}{\mathcal{S}^{\frac{q}{2}}}\Big)^{\frac{2}{4-q}}\Big(\frac{(q-2)\mathcal{S}^2}{b\mathcal{S}^2-\delta}\Big)^{\frac{q-2}{4-q}}
\Big)
\geq 0,
\endaligned
$$
this gives that $f'(t)\geq 0$ for all $t\geq 0$.  Clearly
$f'(t)>0$ for $t>0$ small enough, thus we obtain that
$0=f(0)<f(1)=0$, a contradiction. \vspace{6pt}

\noindent $(ii)$\quad We also use
\begin{equation}\label{eq70}
\overline{I}_0(u)=\frac{a}{2}\|u\|^2-\frac{\lambda}{q}\|u\|_q^q, \ u\in
H^1_0(\Omega)
\end{equation}
to denote the limited functional of $\overline{I}_b$ as $b\to 0^+$
in the following proof.  Clearly, $\overline{I}_0\in C^2(H_0^1(\Omega),\mathbb{R})$ and critical points of
$\overline{I}_0$ are weak solutions of the following equation
$$
\left\{\aligned-a\Delta u &= \lambda|u|^{q-2}u, &\quad \text{in}\ \Omega, \\
u&=0,& \text{on}\ \partial\Omega.
\endaligned
\right.
$$

Noting that if $0<\delta<b\mathcal{S}^2$ then for all $u\in
H^1_0(\Omega)\setminus\{0\}$  we have
$$\aligned
\overline{I}_b(u)=&\frac{a}{2}\|u\|^2+\frac{1}{4}(b\|u\|^4-\delta\|u\|_4^4)-\frac{\lambda}{q}\|u\|_q^q\\
\geq&\frac{a}{2}\|u\|^2+\frac{1}{4}(b-\frac{\delta}{\mathcal{S}^2})\|u\|^4-\frac{\lambda}{q}\|u\|_q^q  \to +\infty
\endaligned
$$
as $\|u\|\to +\infty$, therefore $\overline{I}_b$ is coercive.  On
the other hand, for given $a>0$ and $\lambda>0$, we can find $v_0\in H^1_0(\Omega)\setminus\{0\}$
such that $\overline{I}_0(v_0)<0$, then there exists
$\overline{b}_1>0$ dependent of $a$ and $\lambda$ such that for each $b\in (0, \overline{b}_1)$,
$\overline{I}_b(v_0)<0$, therefore,
$\overline{d}:=\inf\{\overline{I}_b(u): u\in H^1_0(\Omega)\}<0$.  By
Ekeland variational principle, there exists a sequence $\{u_n\}$
satisfying $\overline{I}_b(u_n)\to \overline{d}<0$ and
$\overline{I}'_b(u_n)\to 0$ strongly in $H^{-1}(\Omega)$.  It
follows from Proposition~\ref{pro01} and \eqref{eq69} that there
exists $\overline{u}_1\in H^1_0(\Omega)\setminus\{0\}$ such that
$u_n\to \overline{u}_1$ in $H^1_0(\Omega)$, thus
$\overline{I}'_b(\overline{u}_1)=0$, so that $\overline{u}_1$ is a
solution of \eqref{eq01}.  To obtain the second solution, we just
remind that the functional $\overline{I}_0$ defined by \eqref{eq70}
has the similar properties $(M1)-(M5)$, then by following the proof
of Theorem~\ref{thm04} and combining with Proposition~\ref{pro01}
and \eqref{eq69}, we can get that there exists $\overline{b}_2>0$ dependent of $a$ and $\lambda$ such that for each $b\in (0, \overline{b}_2)$, Eq.~\eqref{eq01} has  a  mountain pass solution $u_2 (\neq u_1)$ in both cases of $\delta<b\mathcal{S}^2$ and $\delta=b\mathcal{S}^2$, here
we omit the details of proof. $\qquad
\raisebox{-0.5mm}{\rule{1.5mm}{4mm}}$\vspace{6pt}

\begin{remark}\label{rmk08}{\em
$\overline{I}_b$ satisfies the $(PS)_c$ condition for all $c\in
\bbr$ if the $(PS)_c$ sequence is bounded, this can be deduced
easily from Proposition~\ref{pro01} and \eqref{eq69}.}
\end{remark}

\noindent{\bf Acknowledgments}

\noindent Part of this work was done while the second author was
visiting the Department of mathematical Sciences, Tsinghua
University, China. He thanks Professor Wen-Ming Zou for the support
and kind hospitality.  Thanks also to Professor K. Perera for his
helpful comments.  This work was supported by Suzhou University of
Science and Technology foundation grant (331412104), Natural Science
Foundation of China (11471235 and 11171247).



\end{CJK*}
\end{document}